\documentclass[12pt]{amsart}
\usepackage{amsmath}
\usepackage{amssymb}
\usepackage{mathrsfs}
\usepackage[a4paper]{geometry}
\theoremstyle{plain}
\newtheorem{teo}{Theorem}
\newtheorem{pro}{Proposition}
\newtheorem{lem}{Lemma}
\newtheorem{cor}{Corollary}

\theoremstyle{definition}

\theoremstyle{remark}

\newtheorem*{nota}{Remark}
\newtheorem*{note}{Remarks}

\newcommand{\p}{\mathbb{P}}
\newcommand{\re}{\mathbb{R}}

\newcommand{\be}{{}_0\mathscr{F}_1}
\newcommand{\fc}{{}_1\mathscr{F}_1}

\newcommand{\jac}{P_n^{r,s}}
\newcommand{\jam}{P_{\tau}^{r,s,\beta}}
\date{\today}
\title{Radial Dunkl Processes : Existence and uniqueness, Hitting time, Beta Processes and Random Matrices} 
\begin{document}
\maketitle
\centerline{NIZAR DEMNI\footnote{Laboratoire de Probabilit\'es et Mod\`eles Al\'eatoires, Universit\'e de Paris VI, 4 Place Jussieu, Case 188, F-75252 Paris Cedex 05,
e-mail : demni@ccr.jussieu.fr}}  

\begin{abstract} 
We begin with the study of some properties of the radial Dunkl process associated to a reduced root system $R$. It is shown that this diffusion is the unique strong solution for all $t \geq 0$ of a SDE with singular drift. Then, we study $T_0$, the first hitting time of the positive Weyl chamber : we prove, via stochastic calculus, a result already obtained by Chybiryakov on the finiteness of $T_0$. The second and new part deals with the law of $T_0$ for which we compute the tail distribution, as well as some insight via stochastic calculus on how root systems are connected with eigenvalues of standard matrix-valued processes. This gives rise to the so-called $\beta$-processes. The ultraspherical $\beta$-Jacobi case still involves a reduced root system while the general case is closely connected to a non reduced one. This process lives in a convex bounded domain known as principal Weyl alcove and the strong uniqueness result remains valid. The last part deals with the first hitting time of the alcove's boundary and the semi group density which enables us to answer some open questions.    
\end{abstract}

\section{Preliminaries}
We begin by pointing out some facts on root systems and radial Dunkl processes. We refer to \cite{Ros1} for the Dunkl theory, to both \cite{Car} and \cite{Hum} for a background on root systems and 
\cite{Chy}, \cite{Gal} for facts on radial Dunkl processes. Let $(V, < , >)$ be a finite real Euclidean space of dimension $m$. A \emph{reduced} root system $R$ is a finite set of non zero vectors spanning $V$ such that : 
\begin{itemize}
\item[1] $R \cap \re \alpha = \{\alpha,-\alpha\}$ for all $\alpha \in R$. \\
\item[2] $\sigma_{\alpha}(R) = R$ \end{itemize} where $\sigma_{\alpha}$ is the reflection with respect to the hyperplane $H_{\alpha}$ orthogonal to $\alpha$ : 
\begin{equation*}
\sigma_{\alpha} (x) =  x - 2\frac{<\alpha, x>}{<\alpha,\alpha>} \alpha , \qquad x \in V \end{equation*}
A simple system $\Delta$ is  a basis of $V$ which induces a total ordering in $R$. A  root $\alpha$ is positive if it is a positive linear combination of elements of $\Delta$. The set of positive roots is called a positive subsystem and is denoted by $R_+$. Note that the choice of $\Delta$ is not unique and that $R_+$ is uniquely determined by $\Delta$.  The reflection group $W$ is the one generated by all the reflections $\sigma_{\alpha}$ for $\alpha \in R$. Recall that $W$ is finite and the only reflections are of the form $\sigma_{\alpha}$ for $\alpha \in R$. 
Given a root system $R$  with associated positive subsytem $R_+$,  let $C$ be the {\it positive Weyl chamber} defined by : 
\begin{equation*}
C := \{x \in V \, < \alpha , x > > 0 \, \forall \, \alpha \in R_+\} = \{x \in V \, < \alpha , x > > 0 \, \forall \, \alpha \in \Delta\} 
 \end{equation*} and $\overline{C}$ its closure. One of the most important properties is that the convex cone $\overline{C}$ is a fundamental domain, that is each $\lambda \in V$ is conjugate to one and only one $\mu \in \overline{C}$. \\
The radial Dunkl process is defined as the $\overline{C}$-valued continuous paths Markov process whose generator is given by : 
\begin{equation*}
\mathscr{L}u(x) = \frac{1}{2}\Delta u(x) + \sum_{\alpha \in R_+}k(\alpha)\frac{< \alpha,\nabla u(x)>}{ < \alpha, x>} \end{equation*} 
with boundary conditions $\nabla u(x)\cdot \alpha = 0$ for all $x \in H_{\alpha},\, \alpha \in R_+$, $k(\alpha) \geq  0 $ is the multiplicity function (invariant under the action of $W$),  
and $u \in C_c^2(C)$. When $k(\alpha) = 1$ for all $\alpha \in R$, we recover the BM constrained to stay in C, studied by Grabiner (\cite{Gra}). 
The semi-group density of $X$ is given by : 
\begin{equation}
\label{sg}
p_t^k(x,y) = \frac{1}{c_kt^{\gamma + m/2}}e^{-(|x|^2+|y|^2)/2t}D_k^W(x,y)\prod_{\alpha \in R_+}<\alpha,y>^{2k(\alpha)}
\end{equation} 
for $x,y \in C$, where $\gamma = \sum_{\alpha\in R_+}k(\alpha)$, 
\begin{equation*}
D_k^W(x,y)  := \sum_{w \in W}D_k\left(\frac{x}{\sqrt t},\frac{wy}{\sqrt t}\right)
\end{equation*}
where $D_k$ denotes the Dunkl kernel and $c_k$ is given by the Macdonald-Mehta integral (\cite{Ros1}). Indeed, as 
$D_k(0,y) = 1$ (\cite{Ros1}), one gets
\begin{equation*}
t^{\gamma + m/2}c_k = |W|\int_{C}e^{-|y|^2/2t}\prod_{\alpha \in R_+}<\alpha,y>^{2k(\alpha)}dy = \int_{\re^m}e^{-|y|^2/2t}\prod_{\alpha \in R_+}|<\alpha,y>|^{2k(\alpha)}dy
\end{equation*}
since $\re^m = \cup_{w \in W}w\overline{C}$. $D_k^W(x,y)$ is known as \emph{the generalized Bessel function} (up to the constant $|W|$). 
This process is obtained by projecting the Dunkl process valued in $\re^m$ (which has right-continuous and left-limits paths, see \cite{Gal}) on $\overline{C}$. 
The latter was already introduced by R\"osler (\cite{Ros1},\cite{Ros2}) and then studied by Gallardo and Yor (\cite{Gal},\cite{Gall}) and Chybiryakov (\cite{Chy}). 

To illustrate all these facts and motivate the reader as well, we will provide some well known examples. We start with the \emph{rank one} case ($m=1$) for which $R = B_1 = \{\pm 1\}$. Hence $k(\alpha) := k_0 \geq 0$ and $X$ is a \emph{Bessel} process (\cite{Rev}) of \emph{index} $\nu = k_0 - 1/2$. When $k_0 > 0$, it is the unique strong solution of : 
\begin{equation*}
dX_t = dB_t + \frac{k_0}{X_t}dt, \quad t \geq 0, \, X_0 = x > 0.
\end{equation*}
where $B$ is a standard BM. Another well known multivariate example is described by the $A_{m-1}$-type root system  defined as :
\begin{equation*}
A_{m-1} = \{\pm (e_i - e_j) \, 1 \leq i < j \leq m\}, \end{equation*} with positive and simple systems given by : 
\begin{equation*}R_+ = \{e_i - e_j, \,1\leq i < j \leq m \} \quad \Delta = \{e_i - e_{i+1}, \,1\leq i  \leq m \}
\end{equation*}
where $(e_i)_{1 \leq i \leq m}$ is the standard basis of $\re^m$. $V$ is the hyperplane of $\re^m$ consisting of vectors that coordinates sum to zero. Without loss of generality, one can take $\re^m$ instead of $V$ so that $C = \{x \in \re^m, x_1 > \dots > x_m\}$. 
Besides, there is only one orbit and $k(\alpha) := k_1 \geq 0$. Thus, the corresponding radial Dunkl process satisfies :  
\begin{equation} \label{Khil}
dX_t^i = d\nu_t^i + k_1 \sum_{j \neq i}\frac{dt}{X_t^i - X_t^j} \qquad 1 \leq i \leq m , \quad t < \tau
\end{equation}
with $X_0^1 > \dots > X_0^m$, where $(\nu^i)_i$ are independent Brownian motions and $\tau$ is \emph{the first collision time} . For strictly positive $k_1$, this process was deeply studied by C\'epa and L\'epingle (\cite{Cep}, \cite{Cepa}, \cite{Cepa1}): it behaves as $m$-interacting particles on the real line with electrostatic repulsions proportional to the inverse of the distance separating them. Moreover, when $k_1=1,1/2$ respectively, this process evolves like eigenvalues process of Hermitian (Dyson model) and symmetric Brownian motions
(\cite{Dyson}, \cite{Gra}). It was shown in \cite{Cep} that this SDE has a unique strong solution for all $t \geq 0$ and $k_1 > 0$. When reading the proof in \cite{Cep}, one hopes to extend this result for any root system since materials used there are not typical for the $A_{m-1}$-type. This was the original motivation of this work. Our first result claims that 
\begin{equation*}
dX_t = dB_t  - \nabla \Phi(X_t) dt , \quad X_0 \in C 
\end{equation*}
where $\Phi(x) = -\sum_{\alpha \in R_+}k(\alpha)\ln(<\alpha,x>), \, k > 0$ , has a unique strong solution for all $t \geq 0$. At the same time and independently, Chybiryakov and Schapira provide two other proofs: both authors used well posed martingale problems associated respectively with the $\re^m$-valued Dunkl and the radial Heckman-Opdam processes as well as geometric arguments (\cite{Chy}, \cite{Scha}).
The curious reader will wonder what happens if $k(\alpha) = 0$ for some $\alpha$? The answer is that the same result holds up to \emph{the first hitting time} of $\partial C$, say $T_0$ (\cite{Chy} p. 37). Next, we are mainly interested in the tail distribution of $T_0$. Before proceeding, we reprove via stochastic calculus that $T_0 < \infty$ if $k(\alpha) < 1/2$ for at least one $\alpha \in R_+$ (see \cite{Chy} for the original proof using local martingales). More precisely, for such an $\alpha$, we prove that almost surely, $<\alpha,X_t> \leq Y_t$  for all $t \geq 0$, where $Y$ is a Bessel process of dimension strictly less than $2$. At this level, other proofs exist for the above results. 
To our knowledge, the contents of the remainder of the paper are new. In \cite{Chy}, the author derived absolute-continuity relations which allow us to write the tail distribution of $T_0$ when starting from $x \in C$. A $W$-invariant analytic $x$-dependent integral, which value at $0$ is given by a Selberg integral, is involved. As far as we know, though $D_k^W(x,y)$ arises as hypergeometric functions for particular root systems (see the end of \cite{Bak}), forward computations are sophisticated and hard. More precisely,  we think that it is possible to use the integral formula given in Corollary 2 in 
\cite{Kan} with the integration range $0 < X_t^1 < \dots < X_t^m < 1$, known as the Macdonald's conjecture, then perform limit and sums operations. The matrix cases for which the Jack parameter equals to $=1,2$ are more handable with the use of properties of zonal polynomials and Schur functions. However, we think that the approach adopted here is more elegant since on one hand, it disgards the special values of the multiplicity function and on the other hand, does not need long hard formulas. It only relies on some properties picked from Dunkl theory. More precisely, it will be shown that  the $x$-dependent integral is an eigenfunction of some operator which involves the generator $\mathscr{L}$ and the so-called \emph{Euler operator} $E_1$. For some particular root systems, this eigenfunction is identified with some hypergeometric series. A surprising fact is that the eigenoperator can be expressed in terms of a Schr\"odinger operator $\mathscr{H}$ and its minimal eigenvalue $E_{min}$ (minimal energy) (see \cite{Ros1} page 18): 
\begin{equation*}
\mathscr{L}  - E_1 := \mathscr{L} - \sum_{i=1}^mx_i\partial_i = -e^{|x|^2/4}(\mathscr{H} - E_{min})e^{-|x|^2/4}\end{equation*}
Moreover, $(X_t)_{t \geq 0}$ specializes for some values of $k$ to eigenvalues processes of self-adjoint matrix processes such as symmetric and Hermitian Brownian motions, Wishart and Laguerre and matrix Jacobi processes. In those cases, computations can be performed using the action of orthogonal and unitary groups. Indeed, Jack polynomials fit zonal polynomials and Schur functions when the Jack parameter equals to $1,2$ respectively (see \cite{Mac}). The two first ones are identified as $A_{m-1}$-type radial Dunkl processes while Wishart and Laguerre processes are related to the $B_m$ root system. The latter goes beyond the radial Dunkl setting: the reduced root system $C_m$ in a particular case (ultraspheric) is involved and more generally, the non reduced system $BC_m$. This connection was deeply investigated in \cite{Bee} while identifying special functions associated with root systems with multivariate hypergeometric series. Among them appear multivariate Gauss hypergeometric series and Jacobi polynomials (\cite{Lass1}) and these are eigenfunctions of the $\beta$-Jacobi generator. The state space is the so-called \emph{principal Weyl alcove} which is now a bounded  convex domain and fundamental for the action of the {\it affine Weyl group}. Hence, the process evolves like particles in an interval. Then, we extend the strong uniqueness Theorem to the Jacobi context. In the remaining part, we derive some properties: we briefly visit  the Brownian motion in the principal Weyl alcove which corresponds to multiplicities all equal to $1$. Then, an analogous result on the finiteness of the first hitting time of alcoves walls is obtained using similar computations as those for $T_0$. Finally, we derive the semi group density and discuss some open questions left in \cite{Dou}.

\section{Radial Dunkl Process : Existence and Uniqueness of a strong solution} 
\begin{teo}
\label{exi}
Let $R$ be a reduced root system. Let:
\begin{equation*}
\Phi(x) = -\sum_{\alpha \in R_+}k(\alpha)\ln(<\alpha,x>) := \sum_{\alpha \in R_+}k(\alpha)\theta(<\alpha,x>), \qquad x \in C 
\end{equation*}
where $k(\alpha) > 0$ for all $\alpha \in R_+$. Then the SDE 
\begin{equation}\label{DE}
dX_t = dB_t - \nabla \Phi(X_t) dt,  \quad X_0 \in C \end{equation}
where $X$ is an adapted continuous process valued in $\overline{C}$ and $B$ is a Brownian motion in $\re^m$, has a unique strong solution. 
\end{teo}
{\it Proof}: From Theorem 2. 2 in \cite{Cepa}, we deduce that : 
\begin{equation}\label{DE1}
dX_t = dB_t - \nabla \Phi(X_t) dt + n(X_t)dL_t, \quad X_0 \in C 
\end{equation} 
where $n(x)$ is a (unitary) inward normal vector to $C$ at $x$ , $L$ is the boundary process defined by: 
\begin{equation*} 
dL_t = {\bf 1}_{\{X_t \in \partial C\}} dL_t, 
\end{equation*} has a unique strong solution for all $t \geq 0$. 
Moreover : 
\begin{eqnarray}\mathbb{E}\left[ \int_0^T {\bf 1}_{\{X_t \in \partial C\}} dt \right] & = & 0 \label{Slim}\\ 
\mathbb{E}\left[\int_0^T |\nabla \Phi(X_t)| dt\right] & < & \infty \label{lotfi}\end{eqnarray}  for all $T > 0$. Thus, it remains to prove that the boundary process vanishes. 
To proceed, we need two Lemmas. 
\begin{nota}
Both Lemmas below discard the reducedness of $R$. In fact, this assumption figures in the definition of the Dunkl process and originates from analytic purposes like the commutativity of Dunkl operators (\cite{Dunkl}).  
\end{nota}

\begin{lem} 
\label{Hamdi} 
Set $dG_t : = n(X_t)dL_t$. Then, $\forall \alpha \in R_+$,  \begin{equation*} {\bf 1}_{\{<X_t,\alpha> = 0\}} < dG_t, \alpha> = 0 \end{equation*}\end{lem}
{\it Proof} : The proof is roughly a generalization of the one in \cite{Cep} for $R= A_{m-1}$ . In order to convince the reader, we provide an outline. 
Using the occupation density formula, we may write ($<\alpha, X> \geq 0$) : 
\begin{equation*}
\int_0^{\infty}L_t^a(<\alpha,X>)\theta^{'}(a)da = <\alpha,\alpha>  \int_0^t \theta^{'}(<\alpha,X_s>) ds \end{equation*}
where $L_t^a(<\alpha,X>)$ is the local time of the real continuous semimartingale $<\alpha, X >$. On the other hand, the following inequaliy holds (instead of (2.5) in \cite{Cep}) for all $a \in C$: 
\begin{align*}
< \nabla \Phi(x),x-a > &= \sum_{\alpha \in R_+}k(\alpha) \theta^{'}(<\alpha,x>)<\alpha,x-a> \\& 
\overset{(1)}{\geq} \sum_{\alpha \in R_+}k(\alpha)[b_{\alpha}\theta^{'}(<\alpha,x>) - c_{\alpha}<\alpha,x-a> - d_{\alpha}] \\&
\geq \min_{\alpha \in R_+}(b_{\alpha}k(\alpha))\sum_{\alpha \in R_+}\theta^{'}(<\alpha,x>) - |x-a|\sum_{\alpha \in R_+}k(\alpha){c_{\alpha}}|\alpha| - 
\sum_{\alpha \in R_+}k(\alpha)d{\alpha}
\\& := A \sum_{\alpha \in R_+}\theta^{'}(<\alpha,x>) - B|x-a| - C
\end{align*}
by Cauchy-Schwarz inequality, where in $(1)$, we used eq. (2.1) in \cite{Cep} : let $g$ be a convex $C^1$-function on an open convex set $D \subset \re^m$, then $\forall a \in D$, there exist
$b,c,d > 0$ such that for all $x \in D$ : 
\begin{equation*}
< \nabla g(x),x-a >  \, \geq \, b|\nabla g(x)| - c |x-a| - d \end{equation*}
 Note also that $A > 0$ since $b_{\alpha}k(\alpha) > 0$ for all $\alpha \in R_+$. Then, the continuity of $X$, (\ref{lotfi}) and the inequality above yield : 
\begin{equation*}
\int_0^t \theta^{'}(<\alpha,X_s>)ds< \infty\end{equation*} which implies that :
\begin{equation*}
\int_0^{\infty}L_t^a(<\alpha,X>)\theta^{'}(a)da < \infty 
\end{equation*}
Thus, $L_t^0(<\alpha,X>) = 0$ since the function $a \mapsto \theta^{'}(a)$ is not integrable at $0$. 
The next step consists in using Tanaka formula to compute  $dZ_t : = <\alpha,X_t> - (<\alpha,X_t>)^+$ for $\alpha \in \Delta$ :
\begin{equation*}
dZ_t = {\bf 1}_{\{<\alpha,X_t> = 0\}}<\alpha,dB_t >  -  {\bf 1}_{\{<\alpha,X_t> = 0\}}< \alpha, \nabla \Phi(X_t) > dt +   {\bf 1}_{\{<\alpha,X_t> = 0\}}<\alpha, dG_t > 
\end{equation*}
It is obvious that the second term vanishes. The first vanishes too since it is a continuous local martingale with null bracket (occupation density formula). 
As $X_t \in \overline{C}$, then $dZ_t = 0$ a.s. which gives the result. $\hfill \blacksquare$

\begin{lem}
\label{Abass}
Let $x \in \partial C$. Then $<n(x), \alpha> \neq 0$ for some $\alpha \in \Delta$ such that $<x,\alpha> = 0$.\end{lem}
{\it Proof} : Let us assume that $<n(x),\alpha> = 0$ for all $\alpha \in \Delta$ such that $<x,\alpha> = 0$. Then, our assumption implies that $<x,\alpha> > 0$ for all $\alpha \in \Delta$ such that $<n(x),\alpha> \neq 0$. If $<n(x),\alpha> < 0$ for these simple roots, then $x- n(x) \in \overline{C}$. By the virtue of the definition of the inward normal $n(x)$ to $C$ at $x$, i. e, 
\begin{equation}
< x - a, n(x) > \, \leq 0, \quad \forall a \in \overline{C}, 
\label{Ines}
\end{equation}  
it follows that $n(x)$ is the null vector which is not possible. Otherwise, choosing : 
\begin{equation*}
0 < \epsilon < \min_{\alpha/<n(x),\alpha >  > 0}\frac{<x,\alpha>}{<n(x),\alpha>}  \end{equation*}
we claim that $a := x - \epsilon n(x) \in \partial C$. Arguing as before, we are done. $\hfill \blacksquare$\\ 
Now we proceed to end the proof of the Theorem. Note first that $\partial C = \cup_{\alpha \in \Delta}H_{\alpha}$ so that : 
\begin{equation*}
 {\bf 1}_{\{X_t \in \partial C\}}dL_t \leq \sum_{\alpha \in \Delta}{\bf 1}_{\{<X_t,\alpha> = 0\}} dL_t. 
\end{equation*}
If $X_t \in H_{\alpha}$ for one and only one $\alpha \in \Delta$. Then, $n(X_t) = \alpha/||\alpha||$ and Lemma \ref{Hamdi} gives
\begin{equation*}   {\bf 1}_{\{<X_t,\alpha> = 0\}} < dG_t, \alpha> = {\bf 1}_{\{<X_t,\alpha> = 0\}} ||\alpha|| dL_t = 0 
\end{equation*}
Hence, the boundary process vanishes. More generally, we can use the inequality above and write 
\begin{align*}
0 \leq L_t & \leq  \sum_{\alpha \in \Delta}\int_0^t{\bf 1}_{\{<X_s,\alpha> = 0\}}{\bf 1}_{\{<n(X_s),\alpha> \neq 0\}} dL_s 
\\& = \sum_{\alpha \in \Delta}\int_0^t{\bf 1}_{\{<n(X_s),\alpha> \neq 0\}} \frac{1}{<n(X_s),\alpha>}{\bf 1}_{\{<X_s,\alpha> = 0\}}<dG_s,\alpha> = 0
\end{align*}
by Lemma \ref{Hamdi}. $\hfill \blacksquare$ 
\begin{nota} 
When $m=1$, $(X_t)_{t \geq 0}$ is a Bessel process of dimension $\delta = 2k_0 + 1$ and  $k_0 > 0 \Leftrightarrow \delta > 1$. It is well known that the local time vanishes 
(see Ch. XI in \cite{Rev}) which fits our result.  
\end{nota}

\section{Finiteness of the first hitting time of the Weyl chamber}
Let $T_0 := \inf\{t > 0, X_t \in \partial C\}$ be the first hitting time of the Weyl chamber. It was shown in \cite{Chy} (see p. 30) that $T_0 = \infty$ almost surely if $k(\alpha) \geq 1/2$ for all 
$\alpha \in R_+$. In \cite{Cepa}, where $R= A_{m-1}$ and $T_0 = \inf\{t > 0,\, X_t^i = X_t^j\, \textrm{for some}\, (i,j)\}$, authors showed that $T_0 < \infty$ a.s. if and only if $0< k_1 < 1/2$. More generally, the following holds (see \cite{Chy} p. 75 for the original proof) :

\begin{pro}\label{Sacha}
Let $\alpha_0 \in \Delta$ and $T_{\alpha_0} := \inf\{t > 0,\, <\alpha_0,X_t> = 0\}$ such that $T_0 = \inf_{\alpha_0 \in \Delta}T_{\alpha_0}$. If $0 < k(\alpha_0) < 1/2$, then  
$(<\alpha_0,X_t>)_{t \geq 0}$ hits almost surely $0$. In particular, $T_0 < T_{\alpha_0} < \infty$ a. s.
 \end{pro}
{\it Proof} : assume $k(\alpha) > 0$ for all $\alpha \in R$ and let $\alpha_0 \in \Delta$. Our scheme is roughly the same as that used in \cite{Cepa}, thus we shall show that the process 
$<\alpha_0,X>$ is almost surely less than or equal to a Bessel process with dimension $2k(\alpha_0) +1$. The result follows from the fact that $2k(\alpha_0) + 1 < 2$ when $k(\alpha) < 1/2$. For this, we use the SDE (\ref{DE}). For all $t \geq 0$,
\begin{align*}
d< \alpha_0, X_t > &= ||\alpha_0|| d\gamma_t + \sum_{\alpha \in R_+}k(\alpha)\frac{<\alpha,\alpha_0>}{<\alpha,X_t >} dt \\&
= ||\alpha_0||d\gamma_t + k_0\frac{||\alpha_0||^2}{< \alpha_0,X_t >}dt + \sum_{\alpha \in R_+\setminus \alpha_0}k(\alpha)\frac{<\alpha,\alpha_0>}{<\alpha,X_t >} dt\end{align*}
where $k_0$ is the value of $k(\alpha_0)$ corresponding to the conjugacy class of $\alpha_0$. Set 
\begin{equation*}
R =  \cup_{j=1}^p R^j  
\end{equation*} where $R^j,\, 1 \leq j \leq p$ denote the conjugacy classes of $R$ under the $W$-action, then 
\begin{equation*}
R_+ = \cup_{i=1}^p R_+^j \end{equation*}
so that: 
\begin{equation*}   
d< \alpha_0, X_t >   = ||\alpha_0||d\gamma_t + k_0\frac{||\alpha_0||^2}{< \alpha_0,X_t >}dt + \sum_{j=1}^p k_j\sum_{\alpha \in R_+^j\setminus \alpha_0}
\frac{<\alpha,\alpha_0>}{<\alpha,X_t >}dt \end{equation*}
For a conjugacy class $R^j$ and $\alpha \in R^j$, if $< \alpha, \alpha_0 > = a(\alpha) > 0$ then, it is easy to check that $< \sigma_0(\alpha), \alpha_0 > = -a(\alpha)$ where $\sigma_0$ is the reflection with respect to the orthogonal hyperplane $H_{\alpha_0}$ defined by : 
\begin{equation*}\sigma_0(x) = x - 2\frac{< x, \alpha_0 >}{<\alpha_0, \alpha_0 >}\alpha_0 
\end{equation*}
Note that $\sigma_0(\alpha)$ belongs to the same conjugacy class of $\alpha$ and that $\sigma_0(\alpha) \in R_+$ for $\alpha \in R_+ \setminus \alpha_0$. 
Indeed, $\sigma_0(R_+\setminus \alpha_0) = R_+\setminus \alpha_0$ for all $\alpha_0 \in \Delta$ (see Proposition 1. 4 in \cite{Hum}). Hence, 
\begin{equation*}
d< \alpha_0, X_t >   = ||\alpha_0||d\gamma_t + k_0\frac{||\alpha_0||^2}{< \alpha_0,X_t >} dt- \sum_{j=1}^p k_j\sum_{\substack{\alpha \in R_+^j\setminus \alpha_0 \\  a(\alpha) >0}}
\frac{a(\alpha)<\alpha - \sigma_0(\alpha), X_t >}{<\alpha,X_t >\, <\sigma_0(\alpha),X_t>}dt 
\end{equation*}
Furthermore, 
\begin{equation*}
\alpha  - \sigma_0(\alpha) = 2\frac{< \alpha,\alpha_0 >}{< \alpha_0, \alpha_0 >}\alpha_0 \quad \Rightarrow \quad < \alpha  - \sigma_0(\alpha) , X_t > 
= 2 a(\alpha)\frac{<\alpha_0, X_t>}{||\alpha_0||^2}   
\end{equation*} 
Consequently, one gets : 
\begin{equation*}
d < \alpha_0, X_t >   = ||\alpha_0||d\gamma_t + k_0\frac{||\alpha_0||^2}{< \alpha_0,X_t >}dt  + F_t \,dt \end{equation*} 
where $F_t < 0$ on $\{T_{\alpha_0} = \infty\}$. Using the comparison Theorem in \cite{Kar} (Proposition 2. 18. p. 293 and Exercice 2. 19. p. 294), one claims that 
$< \alpha_0,X_t >  \, \leq \, Y_{||\alpha_0||^2t}^x$ for all $t \geq 0$ on $\{T_{\alpha_0} = \infty\}$, where $Y^x$ is a Bessel process defined on the same probability space with respect to the same Brownian motion, of dimension $2k_0+1$ and starting at $Y_0 = x \geq \, <\alpha_0,X_0> > 0$. This is not possible since a Bessel process of dimension $< 2$ hits $0$ a. s.  (\cite{Rev} Chap. XI)  $\hfill \blacksquare$
\begin{nota}
If we remove the assumption $k(\alpha) > 0$ for all $\alpha \in R$, then the SDE (\ref{DE}) can be solved up to time $T_0$ when starting from $x \in C$ (see \cite{Chy}). 
\end{nota}
\section{The law of $T_0$}
Here, we focus on the tail distribution of $T_0$ deduced from absolute continuity relations derived in (\cite{Chy}). Recall that (see \cite{Chy}) the \emph{index} of $X$ is defined by 
$l(\alpha) : = k(\alpha) - 1/2$. The last result asserts that if $-1/2 < l(\alpha) < 0$ for some $\alpha \in \Delta$, then $T_0 < \infty$ almost surely. Besides, if $l(\alpha) \geq 0$ for all $\alpha \in \Delta$ then $T_0 = \infty$ almost surely. Taking into account these statements, two major parts are considered: $l(\alpha) \geq 0$ for all $\alpha \in R$ so that the process with index 
$-l$ hits $0$ almost surely, and $l(\alpha) < 0$ for at least one $\alpha$.  
The tail distribution involves a $W$-invariant $x$-dependent integral. Our line of thinking relies on showing that it is an eigenfunction of an appropriate differential operator. Then, using uniqueness results for some differential equations, the tail distribution is written in $A_{m-1}$ and $B_m$ cases by means of multivariate hypergeometric functions. In the last case, we recover known results from matrix theory for Wishart and Laguerre processes. However, we find it better to postpone this in the next section where links with eigenvalues of matrix-valued  processes are detailed.
 
\subsection{A first formula}
Let us denote by $P_x^{l}$ the law of $(X_t)_{t \geq 0}$ starting from $x \in C$.
Let $E_x^{l}$ be the corresponding expectation. Recall that (\cite{Chy}, Proposition 2.15.c), if $l(\alpha) \geq 0 \, \forall \alpha \in R_+$, then: 
\begin{align*}
P_x^{-l}(T_0 > t) &= E_x^l\left[\left(\prod_{\alpha \in R_+}\frac{<\alpha,X_t>}{<\alpha,x>}\right)^{-2l(\alpha)}\right] 
\\& = \prod_{\alpha \in R_+}<\alpha,x>^{2l(\alpha)}\frac{e^{-|x|^2/2t}}{c_kt^{\gamma + m/2}}\int_{C}e^{-|y|^2/2t}D_k^W(\frac{x}{\sqrt t}, \frac{y}{\sqrt t})\prod_{\alpha \in R_+}<\alpha,y> dy
\\& = \prod_{\alpha \in R_+}<\alpha,x>^{2l(\alpha)}\frac{e^{-|x|^2/2t}}{c_kt^{\gamma'}}\int_{C}e^{-|y|^2/2}D_k^W(\frac{x}{\sqrt t}, y)\prod_{\alpha \in R_+}<\alpha,y> dy
\\& : =  \prod_{\alpha \in R_+}<\alpha,x>^{2l(\alpha)}\frac{e^{-|x|^2/2t}}{c_kt^{\gamma'}} g\left(\frac{x}{\sqrt t}\right)
\end{align*}
where $\gamma = \sum_{\alpha \in R_+}k(\alpha)$  and $\gamma' = \gamma - |R_+|/2$. 

Though $D_k^W$ is given by hypergeometric functions in the special cases $A_{m-1}$ and $B_m$ (see the end of \cite{Bak}), the \emph{Jack polynomials} defining them prevent us from making computations. However, this may be possible when these polynomials fit, for some values of $k$, \emph{Zonal polynomials} and \emph{Schur functions} (see \cite{Mac} for definitions). Our main result does not make these restrictions and uses some properties of the Dunkl kernel $D_k$: 

\begin{teo}\label{T1}
Let $T_i$ be the i-th difference Dunkl operator and $\Delta_k = \sum_{i=1}^m T_i^2$ the Dunkl Laplacian (\cite{Ros1}). Define : 
\begin{equation*}
\mathscr{J}_k^x :=  -\Delta_k^x + \sum_{i =1}^m x_i\partial_i^x : = -\Delta_k^x + E_1^x
\end{equation*}
where $E_1^x := \sum_{i=1}^mx_i\partial_i^x$ is the Euler operator and the superscript indicates the derivative action. Then 
\begin{equation*}
\mathscr{J}_k^x \left[e^{-|y|^2/2}D_k^W(x, y)\right] = E_1^y\left[e^{-|y|^2/2}D_k^W(x, y)\right] 
\end{equation*}  
\end{teo}

{\it Proof} : Recall that if $f$ is $W$-invariant then $T_i^xf = \partial_i^xf$ and that $T_i^xD_k(x,y) = y_iD_k(x,y)$ (see \cite{Ros1}). Then, on one hand : 
\begin{align*}
\Delta_k^x D_k^W(x, y) & = \sum_{w \in W}\sum_{i=1}^m(wy)_iT_i^xD_k(x,wy)
 = \sum_{w \in W}\sum_{i=1}^m(wy)_i^2D_k(x,wy) 
 \\& = \sum_{i=1}^my_i^2\sum_{w \in W}D_k(x,wy) : = p_2(y)D_k^W(x,y) 
\end{align*}
On the other hand :	
\begin{align*}
E_1^x D_k^W(x, y) & = \sum_{w \in W}\sum_{i=1}^mx_i T_i^xD_k(x,wy) = \sum_{w \in W}\sum_{i=1}^m(x_i)(wy)_i D_k(x,wy) 
\\& = \sum_{w \in W}<x,wy> D_k(x,wy) = \sum_{w \in W}<w^{-1}x,y> D_k(x,wy) 
\\& = E_1^y D_k^W(x,y)
\end{align*}
where the last equality follows from  $D_k(x,wy) = D_k(w^{-1}x,y)$ since $D_k(wx,wy) = D_k(x,y)$ for all $w \in W$. The result follows from an easy computation.

\begin{cor}
$g$ is an eigenfunction of $-\mathscr{J}_k$ corresponding to the eigenvalue $m+|R_+|$.
\end{cor}
{\it Proof} : Theorem 1 and an integration by parts give : 
\begin{align*}
-\mathscr{J}_k^x g(x) &= - \int_{C}E_1^y\left[e^{-|y|^2/2}D_k^W(x, y)\right]\prod_{\alpha \in R_+}<\alpha,y> dy
\\& = -\sum_{i=1}^m \int_{C}y_i\prod_{\alpha \in R_+}<\alpha,y> \partial_i^y\left[e^{-|y|^2/2}D_k^W(x, y)\right] \,dy
\\& = \sum_{i=1}^m \int_C e^{-|y|^2/2}D_k^W(x, y) \partial_i\left[y_i\prod_{\alpha \in R_+}<\alpha,y>\right] dy
\\& = \int_C  e^{-|y|^2/2}D_k^W(x, y)\prod_{\alpha \in R_+}<\alpha,y> \sum_{i=1}^m \left[1 + \sum_{\alpha \in R_+}\frac{\alpha_iy_i}{<\alpha,y>}\right] dy 
\end{align*}
and the proof ends by summing over $i$. $\hfill \blacksquare$

\begin{itemize}
\item The $A_{m-1}$  case : as mentioned in the introduction, the $A_{m-1}$-type root system is characterized by : 
\begin{eqnarray*}
R  =  \{\pm (e_i - e_j),\, 1 \leq i < j \leq m\}   & &  R_+  =  \{e_i - e_j,\, 1 \leq i < j \leq m\}\\
\Delta  =  \{e_i - e_{i+1}, \, 1 \leq i \leq m\} & & C = \{y \in \re^m, \, y_1 > y_2 > \dots > y_m\}
\end{eqnarray*}
$W = S_m$ is the permutations group and there is one conjugacy class so that $k = k_1 > 0 \Rightarrow \gamma = k_1m(m-1)/2$. 
Moreover, the generalized Bessel function \footnote{Authors use the change of variable $x\mapsto \sqrt{2}x,\, y\mapsto \sqrt{2} y$ to fit the hypergeometric function obtained when deriving the generating function for Hermite polynomials. This in turn will modify the eigenoperator by a multiplying constant (see p. 183).}  is given by (\cite{Bak} p. 212-214,\,\cite{Che}):
\begin{equation*}
\frac{1}{|W|}D_k^W(x,y) =  {}_0F_0^{(1/k_1)}(x,y) := \sum_{p=0}^{\infty}\sum_{\tau}\frac{J_{\tau}^{(1/k_1)}(x)J_{\tau}^{(1/k_1)}(y)}{J_{\tau}^{(1/k_1)}(1)p!}
\end{equation*}
 where $\tau = (\tau_1,\dots, \tau_m)$ is a partition of weight $|\tau| = p$ and length $m$, $J_{\tau}^{(1/k_1)}$ is the Jack polynomial of Jack parameter $1/k_1$\footnote{With the same notations in \cite{Bak}, $k_1 = 2/\alpha$. This can be seen either from the eigenoperator below or from the orthogonality
weight function involved in the semi group density.},
(see \cite{Bak}, \cite{Mac}). Hence, letting $V$ to be the Vandermonde function : 
\begin{align*}
P_x^{-l}(T_0 > t)  &=  V(x)^{2k_0-1}\frac{|W|e^{-|x|^2/2t}}{c_kt^{k_0m(m-1)/2}} \int_{C}e^{-|y|^2/2} {}_0F_0^{(1/k_0)}(\frac{x}{\sqrt t}, y) V(y) dy 
\end{align*} 
Besides, $\mathscr{J}_k$ writes on $W$-invariant functions   
\begin{equation*}
-\mathscr{J}_k^x = D_0^x - E_1^x := \sum_{i=1}^m \partial_i^{2,x} + 2k_1\sum_{i\neq j}\frac{1}{x_i - x_j}\partial_i -  \sum_{i=1}^mx_i \partial_i^x
\end{equation*}
Finally, since $g$ is $W$-invariant, then
\begin{eqnarray*}
\left(D_0^x - E_1^x\right)g(x) & = & m\frac{m+1}{2}g(x), \\
g(0) = \int_C e^{-|y|^2/2} V(y) dy & = & \frac{1}{m!} \int_{\re^m} e^{-|y|^2/2} |V(y)| dy
\end{eqnarray*}
Let us recall that the Gauss hypergeometric function 
\begin{equation*}
{}_2F_1^{(1/k_1)}(e,b,c,z) =  \sum_{p=0}^{\infty}\sum_{\tau}\frac{(e)_{\tau}(b)_{\tau}}{(c)_{\tau}} \frac{J_{\tau}^{(1/k_1)}(z)}{p!}
\end{equation*} 
is the unique symmetric eigenfunction that equals to $1$ at $0$ of (see \cite{Bee} p. 585)  
\begin{equation}\label{GF}
\sum_{i=1}^mz_i(1-z_i)\partial_i^{2,z} + 2k_1\sum_{i\neq j}\frac{z_i(1-z_i)}{z_i - z_j} \partial_i^z+ \sum_{i=1}^m\left[c- k_1(m-1) -\left(e+b+1 - k_1(m-1)\right) z_i\right]\partial_i^z 
\end{equation}
associated to the eigenvalue $meb$. Letting $z = (1/2)(1- x/\sqrt{b}),\,e=(m+1)/2$ and 
\begin{equation*}
c = k_1(m-1) + \frac{1}{2}[e+b+1 - k_1(m-1)] = \frac{b}{2} + \frac{k_1}{2}(m-1) + \frac{m+3}{4} 
\end{equation*}     
the resulting function is an eigenfunction of 
\begin{equation*}
\sum_{i=1}^m (1-\frac{x_i^2}{b})\partial_i^{2,x} + 2k_1\sum_{i\neq j}\frac{(1-x_i^2/b)}{x_i - x_j} \partial_i^x - \sum_{i=1}^m (b+\frac{m+3}{2} - k_1(m-1))\frac{x_i}{b}\partial_i^x 
\end{equation*}
and $D_0^x - E_1^x$ is the limiting operator as $b$ tends to infinity. Hence, 
\begin{pro}For $k_1 \geq 1/2$,  
\begin{equation*}
g(x) = g(0) C(m,k_1)\lim_{b \rightarrow \infty} {}_2F_1^{(1/k_1)}\left[m+1, b, \frac{b}{2} + \frac{k_1}{2}(m-1) + \frac{m+3}{2}, \frac{1}{2}\left(1-\frac{x}{\sqrt{b}}\right)\right]
\end{equation*}
where 
\begin{equation*}
C(m, k_1)^{-1} = \lim_{b \rightarrow \infty}{}_2F_1^{(1/k_1)}\left(m+1, b, \frac{b}{2} + \frac{k_1}{2}(m-1) + \frac{m+3}{2}, \frac{1}{2}\right)
\end{equation*}\end{pro}

\item The $B_m$ case : This root system is defined by 
\begin{eqnarray*}
R  =  \{\pm e_i, \pm e_i \pm e_j,\, 1 \leq i < j \leq m\}   & &  R_+  =  \{e_i, 1\leq i \leq m, \, e_i \pm e_j,\, 1 \leq i < j \leq m\}\\
\Delta  =  \{e_i - e_{i+1}, \, 1 \leq i \leq m, \, e_m \} & & C = \{y \in \re^m, \, y_1 > y_2 > \dots > y_m > 0\}
\end{eqnarray*}
\noindent 
The Weyl group is generated by transpositions and  '' change sign'' reflections $(x_i \mapsto -x_i)$ and there are two conjugacy classes so that $k=(k_0,k_1) \Rightarrow 
\gamma = mk_0 + m(m-1)k_1$. The generalized Bessel function\footnote{there is an erroneous sign in one of the arguments in \cite{Bak}. Moreover, to recover this expression in the $B_m$ case from that given in \cite{Bak}, one should make substitutions $a = k_0 - 1/2, \, k_1 = 1/\alpha,\,  q = 1 +(m-1)k_1$. We point to the reader that this is different from the one used in \cite{Che} p. 121. } 
is given by (\cite{Bak} p. 214) :  
\begin{equation*}
\frac{1}{|W|}D_k^W(x,y) = {}_0F_1^{(1/k_1)}(k_0 + (m-1)k_1 + \frac{1}{2}, \frac{x^2}{2t},\frac{y^2}{2t})
\end{equation*} 
where  
\begin{equation*}
{}_0F_1^{(1/k_1)}(c,x,y) =  \sum_{p=0}^{\infty}\sum_{\tau}(c)_{\tau}\frac{J_{\tau}^{(k_1)}(x)J_{\tau}^{(1/k_1)}(y)}{J_{\tau}^{(1/k_1)}(1)p!}
\end{equation*}
and $(c)_{\tau} : = \prod_{i=1}^m(c- k_1(i-1))_{\tau_i}$ is the generalized Pochammer symbol (see \cite{Bak}). Then, one has : 
\begin{equation*}
g(x) = |W|\int_C e^{-|y|^2/2}{}_0F_1^{(1/k_1)}(k_0 + (m-1)k_1 + \frac{1}{2}, \frac{x^2}{2},\frac{y^2}{2})\prod_{i=1}^m(y_i)V(y^2)dy
\end{equation*}
The eigenoperator writes on $W$-invariant functions: 
\begin{eqnarray*}
-\mathscr{J}_k^x &=& \sum_{i=1}^m\partial_i^{2,x} + 2k_0 \sum_{i=1}^m\frac{1}{x_i}\partial_i^x + 2k_1\sum_{i \neq j}\left[\frac{1}{x_i - x_j} + \frac{1}{x_i+x_j}\right] \partial_i^x - E_1^x\\
-\mathscr{J}_k^x g(x) & = & m(m+1) g(x),\, \quad g(0) = \frac{1}{2^m m!}\int_{\re^m}e^{-|y|^2}\prod_{i=1}^m|y_i| |V(y^2)| dy.
\end{eqnarray*}
 A change of variable $x_i = \sqrt{2y_i}$ shows that $u(y) : = g(\sqrt{2y})$ satisfies 
 \begin{eqnarray*}
-\tilde{\mathscr{J}}_k^y u(y)  & = & m\frac{(m+1)}{2} u(y),\, \quad g(0) = u(0)\\
-\tilde{\mathscr{J}}_k^y  & = & \sum_{i=1}^m y_i\partial_i^{2,y} + 2k_1\sum_{i \neq j} \frac{y_i}{y_i - y_j}\partial_i^y + \left(k_0 + \frac{1}{2}\right)\sum_{i=1}^m\partial_i^y - E_1^y.
\end{eqnarray*} 
which implies that : 
\begin{equation*}
u(y) = u(0){}_1F_1^{(1/k_1)}(\frac{m+1}{2}, k_0 + (m-1)k_1 + \frac{1}{2}, y)
\end{equation*}
where 
\begin{equation*}{}_1F_1^{(1/k_1)}(b,c,z) =  \sum_{p=0}^{\infty}\sum_{\tau}\frac{(b)_{\tau}}{(c)_{\tau}} \frac{J_{\tau}^{(1/k_1)}(z)}{p!}\end{equation*}
This can be seen from the differential equation (\ref{GF}) and using (\cite{Bak}): 
\begin{equation*}
\lim_{e \rightarrow \infty}{}_2F_1^{(1/k_1)}(e,b,c,\frac{z}{e}) = {}_1F_1^{(1/k_1)}(b,c,z)
\end{equation*}
Finally 
\begin{equation*}
g\left(\frac{x}{\sqrt t}\right) = g(0){}_1F_1^{(1/k_1)}(\frac{m+1}{2}, k_0 + (m-1)k_1 + \frac{1}{2}, \frac{x^2}{2t})
\end{equation*}

Hence, the tail distribution is given by : 
\begin{pro}\label{B1} For $k_0,k_1 \geq 1/2$, 
\begin{equation*}
P_x^{-l}(T_0 > t) = C_k\prod_{i=1}^m \left(\frac{x_i^2}{2t}\right)^{k_0-1/2}\left(V\left(\frac{x^2}{2t}\right)\right)^{2k_1-1}e^{-|x|^2/2t}
{}_1F_1^{(1/k_1)}(\frac{m+1}{2}, k_0 + (m-1)k_1 + \frac{1}{2} , \frac{x^2}{2t}) 
\end{equation*}\end{pro}

\begin{nota} 
1/Adopting the notations used in \cite{Bak}, one has : 
\begin{eqnarray*}
-\tilde{\mathscr{J}}_k^y &=& D_1^y + (a+1)E_0^y - E_1^y \quad (R=B_m, \, y = x^2),
\end{eqnarray*}
Besides, Theorem \ref{T1} was derived there differently for both $A_{m-1}$ and $B_m$ cases when proving a generating function Theorem for generalized Hermite and Laguerre polynomials (page 183 and 192, see also \cite{Che}).  
\end{nota}
\end{itemize}

\subsection{A second formula}
In \cite{Chy} (see Proposition 2.15.b), the author derived another absolute-continuity relation from which we deduce that if $l(\alpha) < 0$ for at least one $\alpha \in R_+$, then 
\begin{align*}
P_x^l( T_0 > t) &= E_x^0\left[\prod_{\alpha \in R_+}\left(\frac{<\alpha,X_t>}{<\alpha,x>}\right)^{\l(\alpha)}\exp\left(-\frac{1}{2}\sum_{\alpha,\gamma \in R_+}
\int_0^t\frac{<\alpha,\gamma>l(\alpha)l(\gamma)}{<\alpha,X_s><\gamma,X_s>}ds\right)\right]
\\& = E_x^r\left[\prod_{\alpha \in R_+}\left(\frac{<\alpha,X_t>}{<\alpha,x>}\right)^{\l(\alpha) - r(\alpha)}\exp\left(-\frac{1}{2}\sum_{\alpha,\gamma \in R_+}
\int_0^t\frac{<\alpha,\gamma>l(\alpha,\gamma)}{<\alpha,X_s><\gamma,X_s>}ds\right)\right]
\end{align*} 
where the last equality follows from part (c) of the same Proposition, $l(\alpha,\gamma) = l(\alpha)l(\gamma) - r(\alpha)r(\gamma)$ and  
\begin{equation*}
r(\alpha) = \left\{
\begin{array}{lcl}
l(\alpha) & \textrm{if} & l(\alpha) \geq 0\\
-l(\alpha) & \textrm{if} & l(\alpha) < 0
\end{array}\right. 
\end{equation*}
Then $l(\alpha,\gamma) = 0$ if $l(\alpha)l(\gamma) \geq 0$ and $l(\alpha,\gamma) = -2r(\alpha)r(\gamma)$ else. As a result, 
\begin{equation*}
P_x^l(T_0 > t) = E_x^r\left[\prod_{\substack{\alpha \in R_+\\l(\alpha) < 0}}\left(\frac{<\alpha,x>}{<\alpha,X_t>}\right)^{ 2r(\alpha)}
\exp\left(\sum_{\substack{\alpha,\gamma \in R_+\\l(\alpha)l(\gamma) < 0}} \int_0^t\frac{<\alpha,\gamma>r(\alpha)r(\gamma)}{<\alpha,X_s><\gamma,X_s>}ds\right)\right]
\end{equation*}
Next, note that the exponential functional equals $1$ for both root systems $A_{m-1}$ and $B_m$. For the first, it is obvious since $R$ consists of one orbit 
so that $\{\alpha,\,\gamma \in R_+, l(\alpha) l(\gamma) < 0\}$ is empty. This gives the same expression already considered in the previous subsection. For the second, writing 
$R_+ = \{e_i, 1 \leq i \leq m\} \cup \{e_j \pm e_k, 1 \leq j < k \leq m\}$ so that 
$<e_i,e_j \pm e_k> = \delta_{ij} \pm \delta_{ik}$ gives  
\begin{align*}
& S =  \sum_{i=1}^m\sum_{i < k}\frac{1}{X_t^i}\left[\frac{1}{X_t^i - X_t^k} + \frac{1}{X_t^i + X_t^k}\right]  + \sum_{i=1}^m\sum_{k < i}\frac{1}{X_t^i}\left[\frac{-1}{X_t^k - X_t^i} 
+ \frac{1}{X_t^k + X_t^i}\right]
\\& = \sum_{i=1}^m\sum_{i < k}\frac{2}{(X_t^i)^2 - (X_t^k)^2} - \sum_{i=1}^m\sum_{k < i}\frac{2}{(X_t^k)^2 - (X_t^i)^2} = 0
\end{align*}
where $S$ stands for the sum between brackets. The reader can also check that this holds for $C_m$ and $D_m$ root systems (see the end of the paper for definitions). 
However, we restrict ourselves to the $B_m$ -case since, for particular values of the multiplicity function, we will recover a known result from matrix theory (see next section). Let us investigate the case $k_0 < 1/2,\,k_1 \geq 1/2$ for which $l_0 < 0, \, l_1 \geq 0$. One writes : 
\begin{align*}
g(x) =  \int_C e^{-|y|^2/2}D_k^W(x,y)\prod_{i=1}^m(y_i)V^{2k_1}(y^2)dy
\end{align*}  
The machinery used before still applies and gives : 
\begin{equation*}
-\mathscr{J}_kg = 2m[ 1+ k_1(m-1)]g
\end{equation*} 
Thus
\begin{pro} 
In the $B_m$ case and for $k_0 <1/2,\, k_1 \geq 1/2$, one has:
\begin{equation*}
P_x^{l}(T_0 > t) = C_k\prod_{i=1}^m \left(\frac{x_i^2}{2t}\right)^{k_0-1/2}e^{-|x|^2/2t}{}_1F_1^{(1/k_1)}(1+k_1(m-1), k_0 + (m-1)k_1 + \frac{1}{2} , \frac{x^2}{2t}) 
\end{equation*}\end{pro}
In the remaining case $k_0 \geq 1/2,\,k_1 < 1/2$, the tail distribution writes : 
\begin{align*}
g(x) =  \int_C e^{-|y|^2/2}D_k^W(x, y)\prod_{i=1}^m(y_i)^{2k_0}V(y^2)dy
\end{align*}   
Thus 
\begin{equation*}
-\mathscr{J}_k g = m[2k_0 + m]g
\end{equation*}
so that
\begin{pro}For $k_0 \geq 1/2,k_1 <1/2$,
\begin{equation*}
P_x^{l}(T_0 > t) = C_k V\left(\frac{x^2}{2t}\right)^{2k_1 - 1}e^{-|x|^2/2t}{}_1F_1^{(1/k_1)}(k_0 + \frac{m}{2}, k_0 + (m-1)k_1 + \frac{1}{2}, \frac{x^2}{2t}) 
\end{equation*}\end{pro}

\section{$\beta$-processes and random matrices}
In the sequel, we will see how eigenvalues of some classical matrix-valued processes and radial Dunkl processes are interelated using SDE. 
This connection was already checked by physicists throughout eigenvalues probability densities and Fokker-Planck equations for parameter-dependent random matrices (\cite{Che}).    
As we mentioned in the introduction, the $A_{m-1}$-type is connected to symmetric and Hermitian Brownain motions. Set $k := \beta/2,\, \beta > 0$, then such a process will be called $\beta$-Dyson, referring to the Dyson model when $\beta = 2$. This parameter is called the Dyson index. Henceforth, we will adopt new notation for the eigenvalues process, we wil write $\lambda$ instead of $X$.

\subsection{The $B_m$-type: $\beta$-Laguerre processes.} 
The $B_m$ system turns out to be related to eigenvalues of Wishart and Laguerre processes which satisfy the following stochastic differential system (see \cite{Bru},\cite{Dem}): 
\begin{equation*}
d\lambda_i(t) = 2\sqrt{\lambda_i(t)}\, d\nu_i(t) + \beta\left[ \delta + \sum_{k \neq i}\frac{\lambda_i(t) + \lambda_k(t)}{\lambda_i(t) - \lambda_k(t)}\right]dt \qquad 1 \leq i \leq m. \end{equation*}for $\beta = 1, 2$ and $\delta \geq m+1, m$ respectively, where $(\nu_i)_i$ are independent Brownian motions and $\lambda_1(0) > \dots > \lambda_m(0)$. Recall that the process remains strictly positive if it is initially strictly positive. This suggests to define the $\beta$-Laguerre process as the solution, when it exists, of : 
\begin{equation*}
d\lambda_i(t) = 2\sqrt{\lambda_i(t)}\, d\nu_i(t) + \beta\left[ \delta + \sum_{k \neq i}\frac{\lambda_i(t) + \lambda_k(t)}{\lambda_i(t) - \lambda_k(t)}\right]dt \qquad 1 \leq i \leq m,
\qquad t < \tau \wedge R_0, 
\end{equation*} 
where $R_0 = \inf\{t , \, \lambda_m(t) = 0\}$,  $\beta,\delta > 0$ and with $\lambda_1(0) > \dots > \lambda_m(0) > 0$. 
It is very easy to see that : 
\begin{equation*}
dD_t := d\left(\prod_{i=1}^m\lambda_i(t)\right) = 2D_t \sqrt{\sum_{i=1}^m\frac{1}{\lambda_i(t)}} d\tilde{B}_t + \beta (\delta-m+1) D_t \sum_{i=1}^m\frac{1}{\lambda_i(t)}dt \end{equation*}
for all $t < \tau \wedge R_0$ where $\tilde{B}$ is a standard Brownian motion. It follows that $\forall \, r \in \re$: 
\begin{eqnarray*}
d(\ln(D_t)) & = & 2 \sqrt{\sum_{i=1}^m\frac{1}{\lambda_i(t)}} d\tilde{B}_t + [\beta (\delta - m) + \beta - 2]  \sum_{i=1}^m\frac{1}{\lambda_i(t)}dt \\
d(\det(D_t)^r) & = & M_t + r[\beta(\delta-m+1) + 2r -2]D_t^r \sum_{i=1}^m\frac{1}{\lambda_i(t)}dt  \end{eqnarray*} 
where $M_t = 2r D_t^r \sqrt{\sum_{i=1}^m{1/\lambda_i(t)}} d\tilde{B}_t $. From these two SDE, we can argue as in the Wishart and Laguerre cases that $R_0 > \tau$ a.s. when $\delta \geq m-1+2/\beta$ (choose $2r = 2- \beta(\delta-m+1) < 0$ when $\delta > m-1 + 2/\beta$ then use McKean's argument). Set $r_i := \sqrt{\lambda_i}$, then, for $t <  \tau \wedge R_0$: 
\begin{align*}dr_i(t) &= d\nu_i(t) + \frac{1}{2r_i(t)}\left[\beta \delta -1 + \beta \sum_{j \neq i}\frac{r_i^2 + r_j^2}{r_i^2-r_j^2}\right] dt \\&
= d\nu_i(t) + \frac{\beta (\delta-m+1)-1}{2r_i(t)}dt + \frac{\beta}{2} \sum_{j \neq i}\left[\frac{1}{r_i(t) - r_j(t)} + \frac{1}{r_i(t) + r_j(t)}\right] dt \\&
= d\nu_i(t) + \frac{k_0}{r_i(t)}dt + k_1 \sum_{j \neq i}\left[\frac{1}{r_i(t) - r_j(t)} + \frac{1}{r_i(t) + r_j(t)}\right] dt
\end{align*}
with $2k_0 = \beta(\delta - m + 1) -1,\, 2k_1 = \beta$.   
Consequently, the process $r = (r_1,\dots,r_m)$ defined for all $t < \tau \wedge R_0$ is a $B_m$-radial Dunkl process. 
Using Theorem \ref{exi}, one claim that the SDE above has a unique strong solution for all $t \geq 0$ and all $\beta, \, \delta$ such that $k_0, \, k_1 > 0$.
This strengthen results from matrix theory : in the Wishart setting ($\beta =1$), the strong uniqueness holds for $\delta > m$ and in the Laguerre case ($\beta = 2$), it holds for 
$\delta > m - 1/2$.    
Besides, the generalized Bessel function is given by (\cite{Bak})\footnote{With the same notations used in \cite{Bak}, one has $\beta a'/2 = k_0,\, a = k_0 - 1/2,\, \beta = 2/\alpha 
\Rightarrow a+q = \beta\delta/2$.} : 
\begin{equation*}
\frac{1}{|W|}\sum_{w \in W}D_k(x,wy) =  {}_0F_1^{(2/\beta)}(\frac{\beta \delta}{2}, \frac{x^2}{2},\frac{y^2}{2}) 
: = \sum_{p=0}^{\infty}\sum_{\tau} \left(\frac{\beta\delta}{2}\right)_{\tau} \frac{J_{\tau}^{(2/\beta)}(x^2/2)J_{\tau}^{(2/\beta)}(y^2/2)}{J_{\tau}^{(2/\beta)}(1_m)p!}
\end{equation*}
so that (\ref{sg}) writes 
\begin{equation}\label{LE}
p_t^{k_0,k_1}(x,y) = \frac{|W|}{c_kt^{\gamma + m/2}}e^{-(|x|^2+|y|^2)/2t}{}_0F_1^{(2/\beta)}(\frac{\beta \delta}{2}, \frac{x^2}{2t},\frac{y^2}{2t})
\prod_{i=1}^m(y_i)^{2k_0}V^{2k_1}(y^2)dy 
\end{equation}
where $V$ stands for the Vandermonde function. Using the variable change $y \mapsto \sqrt y$,  the semi-group density of the $\beta$-Laguerre process writes : 
\begin{equation*}
q_t^{k_0,k_1}(x,y) = \frac{C_{k_0,k_1}}{t^{\gamma +2k_1+m/2}} e^{-\left(\sum_{i=1}^m(x_i+y_i)/2t\right)} {}_0F_1^{(2/\beta)}(\frac{\beta \delta}{2}, \frac{x}{2t},\frac{y}{2t})
\prod_{i=1}^m(y_i)^{k_0-1/2}V^{2k_1}(y)dy 
\end{equation*}
For $x=0$ and $t=1$, we recover the same p.d.f. given in \cite{Dum} for $\beta$-Laguerre ensemble.

\begin{note} 
1/ Recall that for all $\alpha \in R$, we set $l(\alpha) = k(\alpha) - 1/2$. Hence, in the $B_m$-case, $l_0 = k_0 - 1/2,\,l_1 = k_1 - 1/2$. For $-l$, all corresponding parameters will be primed. For instance, $-l_0 = k'_0 - 1/2, \, -l_1= k'_1 - 1/2$.  Let us consider a  Wishart process of dimension $\delta'$ such that $m-1 \leq \delta' < m+1$ (\cite{Bru}), 
$k_1' = 1/2 \, (\beta' =1)$ and $k_0' = (\delta' - m)/2 \Rightarrow -l_1 = 0,\, -l_0 = (\delta' -m - 1)/2 < 0$. 
Set $\delta' = m+1 -2\nu,\, 0 < \nu < 1/2$, then, $l_1 = 0, l_1 = \nu \Rightarrow k_1 = 1/2\,(\beta = 1)$ and $k_0 = \nu + 1/2 \, (\delta = m+1+2\nu)$. Results of 4.1 writes: 
\begin{equation*}
P_x^{-l}(T_0 > t) = C_k\prod_{i=1}^m \left(\frac{x_i^2}{2t}\right)^{\nu}e^{-|x|^2/2t}{}_1F_1^{(2)}(\frac{m+1}{2}, \frac{\delta}{2}, \frac{x^2}{2t}) 
\end{equation*} 
which fits the expression already derived in \cite{Don}. When $k'_0 = k'_1 = 0$ ($-l_0 = -l_1 = -1/2$), then $k_0 = k_1 = 1 (\beta = 2, \delta = m+1/2)$ and the Jack polynomials fits the Schur functions (see \cite{Mac}). In that case, the following representation holds (\cite{Gro})  
\begin{equation*}
{}_1F_1^{(1)}(a, b, z) = \frac{\det(z_i^{m-j}\fc(a - j +1, b-j+1, z_i)_{1 \leq i,j \leq m}}{V(z)} 
\end{equation*}
where $\fc$ denotes the univariate hypergeometric function. Hence, the tail distribution writes :
\begin{equation*}
P_x^{-l}(T_0 > t) = C_k\det\left[\left(\frac{x_i^2}{2t}\right)^{m-j+1/2}e^{-x_i^2/2t}\fc(\frac{m+1}{2}-j+1,m-j+\frac{3}{2},\frac{x_i^2}{2t})\right]_{1\leq i,j, \leq m} 
\end{equation*}
The corresponding process is the Brownian motion in the Weyl chamber of $B$-type.\\
For Laguerre processes (\cite{Dem}) of dimension $\delta$ such that $m-1/2 \leq \delta < m$, one should apply results derived in section 4.2. 
Take $k_1' = 1\, (\beta' = 2)$ and $k_0' = \delta'-m + 1/2 \Rightarrow l_1 = 1/2,\, l_0 = \delta' - m := -\nu$ with $0< \nu < 1/2$. Thus $r_1 = 1/2,\, r_0 = \nu \Rightarrow k_1 = 1 \,(\beta =2)$ and $k_0 = \nu +1/2 \, (\delta = m+\nu)$ so that : 
\begin{equation*}
P_x^{l}(T_0 > t) = C_k\prod_{i=1}^m \left(\frac{x_i^2}{2t}\right)^{\nu}e^{-|x|^2/2t} {}_1F_1^{(1)}(m, \delta, \frac{x^2}{2t}) 
\end{equation*}

2/Recall that when $\beta=2$, ${}_0F_1^{(1)}$ has a determinantal representation (see \cite{Gro}) yielding to K\"onig and O'Connell result on the $V$-transform of $m$-independent squared Bessel processes (BESQs) constrained never to collide (or stay in the $A_{m-1}$-type Weyl chamber, see \cite{Konig}). Similar results holds for $A_{m-1}$-type root system with Brownian motions instead of BESQs. Nevertheless, when $k_0 = k_1 = 1\, (\beta = 2, \delta = m+1/2)$, a similar interpretation involving $m$- independent Brownian motions killed when they reaches $0$ holds. However, the Vandermonde function may be replaced by the product over positive roots. In this case, the eigenvalues process is known as the BM in the Weyl chamber of type $B_m$ (see \cite{Gra} for further details and other root systems). Since $\gamma = m^2$ and from (\cite{Gro}, \cite{Dem}) : 
\begin{equation*}
{}_rF_s^{(1)}((m+a_i)_{1\leq i \leq r}, (m+b_j)_{1 \leq j \leq s}, x,y) = \frac{\det[{}_r\mathscr{F}_s((a_i+1)_{1 \leq i \leq r},(b_j+1)_{1 \leq j \leq s}, x_ly_f)]_{l,f}}{V(x)V(y)}
\end{equation*}
$(\ref{LE})$ transforms to : 
\begin{align*}
p_t^{1,1}(x,y) & = C_m\frac{h(y)}{h(x)}\frac{e^{-(|x|^2+|y|^2)/2t}}{t^{m/2}}\prod_{i,j=1}^m\left(\frac{x_iy_j}{t}\right)\det\left[\be\left(\frac{1}{2}+1, \frac{(x_iy_j)^2}{4t^2}\right)\right]_{i,j}
\\& = C_m \frac{h(y)}{h(x)}\frac{e^{-(|x|^2+|y|^2)/2t}}{t^{m/2}}\det\left[\frac{x_iy_j}{t}\be\left(\frac{3}{2}, \frac{(x_iy_j)^2}{4t^2}\right)\right]_{i,j}
\end{align*}
where $h$ is the product over positive roots. Besides, the following holds (see \cite{Bry3}) : 
\begin{equation*}
\be(\frac{3}{2},z) =  \frac{C}{2\sqrt{z}}\sinh(2\sqrt{z}).
\end{equation*}    
Thus, 
\begin{align*}
p_t^{1,1}(x,y) & = C_m\frac{h(y)}{h(x)}\frac{1}{(2\pi t)^{m/2}}e^{-(|x|^2+|y|^2)/2t}\det\left[\sinh\left(\frac{x_iy_j}{t}\right)\right]_{i,j}
\\& = \frac{h(y)}{h(x)}\det\left[N_t(y_j - x_i) - N_t(y_j+x_i)\right]_{i,j}
\end{align*}
where $N_t(u) = (1/\sqrt{2\pi t})e^{-u^2/2t}$, which fits Grabiner's result (\cite{Gra} page 186). This is in agreement with the generator since $\Delta h = 0$ (\cite{Gra}) and 
\begin{equation*}
\mathscr{L}f = \Delta f + \Gamma(\log h, f) = \Delta f + \sum_{i=1}^m\partial_i(\log h)\partial_i f,
\end{equation*}
where $\Gamma$ is the so-called ''op\'erateur du carr\'e du champ'' (see \cite{Rev} Chap. VIII). Besides, for $m=1$, $r$ is a Bessel process of dimension $2\delta = 3$ and the expression inside the determinant in the second line is exactly the semi-group of the Brownian motion killed when it reaches $0$ (see \cite{Rev} p. 87).\\ 
\end{note}

\subsection{Generalized Bessel function in the $D_m$ case}
In the classification of root systems, the $A_{m-1}$ and $B_m$ are known to be ''irreducible'' and both of them correspond to some matrix processes. Another one, yet with no underlying matrices, is the $D_m$ root system defined by (see \cite{Hum} p. 42)  
\begin{equation*}
R  = \{\pm e_i \pm e_j, \,1 \leq i < j \leq m\},\qquad R_+ = \{e_i \pm e_j, \,1 \leq i < j \leq m\}
\end{equation*}
There is one conjugacy class so that $k(\alpha) = k_1$. Grabiner's result reads for the Brownian motion in the Weyl chamber of $D_m$-type ($k_1 = 1$) : 
\begin{align*}
p_t^1(x,y) &= \frac{V(y^2)}{V(x^2)}\frac{\det[N_t(y_i - x_j) - N_t(y_i + x_j)] + \det[N_t(y_i - x_j) + N_t(y_i + x_j)]}{2}
\\ & = \frac{C_m}{t^{\gamma + m/2}}e^{-(|x|^2+|y|^2)/2t}\frac{\det\left[\sinh(x_iy_j/t)\right] + \det\left[\cosh(x_iy_j/t)\right]}{V(x^2/4t^2)V(y^2)}V^2(y^2)
\end{align*}   
where $\gamma = m(m-1)$. The second term in the sum involves the transition density of a reflected Brownian motion ($|B|$, see \cite{Rev} p. 81). A natural way to interpret the announced formula is that the Weyl chamber is given by : 
\begin{equation*}
C = \{x \in \re^m,\, x_1 > \dots > x_{m-1} > |x_m|\},
\end{equation*} 
so that $C$ fits the $B_m$-Weyl chamber when $x_m > 0$, otherwise, it is its conjuguate with respect to $s_{e_m}$ since this simple reflection acts only on $x_m$ and retains the others. With the help of the determinantal formula used before (\cite{Gro}, \cite{Dem}), $\be(\frac{3}{2},z) =  C\sinh(2\sqrt{z})/\sqrt{z}$ and $\be(1/2, \, z) = \cosh(2\sqrt{z})$ (\cite{Bry3}), one writes:  
\begin{align*}
p_t^1(x,y) =  \frac{e^{-(|x|^2+|y|^2)/2t}}{c_k t^{\gamma + m/2}}
\left[\prod_{i=1}^m \left(\frac{x_iy_i}{2t}\right){}_0F_1^{(1)}\left(m+\frac{1}{2}, \frac{x^2}{2t},\frac{y^2}{2t}\right)+ {}_0F_1^{(1)}\left(m-\frac{1}{2}, \frac{x^2}{2t},\frac{y^2}{2t}\right)\right] 
V^2(y^2)\end{align*}
With regard to (\ref{sg}) and setting $q = 1+(m-1)k_1$, it is natural to prove that : 
\begin{pro}
\begin{equation*}
\frac{1}{|W|}\sum_{w \in W} D_k(x,wy) = \prod_{i=1}^m \left(\frac{x_iy_i}{2}\right){}_0F_1^{(1/k_1)}\left(q+\frac{1}{2}, \frac{x^2}{2},\frac{y^2}{2}\right)+ 
{}_0F_1^{(1/k_1)}\left(q-\frac{1}{2}, \frac{x^2}{2t},\frac{y^2}{2t}\right).
\end{equation*}
\end{pro}
{\it Proof}: it relies on the following characterization (\cite{Ros1}): given a reduced root system $R$ with finite reflection group $W$, $D_k^W(x,\cdot)$ is the unique $W$-invariant function valued $1$ at $x=0$ satisfying 
$\Delta_kD_{k}^W(x,\cdot) = |x|^2D_k^W(x,\cdot)$. 
It is easy to see that the function above is $W$-invariant since $W$ is the semi-direct product of the symmetric group $S_m$ and $(\mathbb{Z}/2\mathbb{Z})^{m-1}$ acting by even sign changes. However, it is not for the finite reflection group associated to the $B_m$ root system due to the term multiplying the first hypergeometric series. In the $D_m$ case, the Dunkl Laplacian writes on $W$-invariant functions : 
\begin{equation*}
\Delta_k = \sum_{i=1}^m\partial_i^2 + 2k_1\sum_{i\neq j}\left[\frac{1}{y_i - y_j} + \frac{1}{y_i+y_j}\right]\partial_i
\end{equation*}
Let : 
\begin{eqnarray*}
f(x,y) &=& {}_0F_1^{(1/k_1)}\left(q+\frac{1}{2}, \frac{x^2}{2},\frac{y^2}{2}\right)\\
g(x,y) &=& {}_0F_1^{(1/k_1)}\left(q-\frac{1}{2}, \frac{x^2}{2},\frac{y^2}{2}\right),\\ 
d(x,y) &=& \prod_{i=1}^m(x_iy_i/2)
\end{eqnarray*} 
considered as functions of the variable $y$ such that the generalized Bessel function is proportional to $df + g$ and $\Delta_k[df + g] = \Delta_k(df) + \Delta_k(g)$. Recall that:  
\begin{equation*}
\frac{1}{|W|}\sum_{w \in W}D_k^{(B_m)}(x,wy) =  {}_0F_1^{(1/k_1)}(k_0 - \frac{1}{2} + q, \frac{x^2}{2},\frac{y^2}{2})
\end{equation*} 
Note also that $\Delta_k$ is a particular case of the Dunkl Laplacian considered for the $B_m$ root system when $k_0= 0$. As a result 
\begin{equation*}
\Delta_k g(x,y) = \Delta_k^{(B_m)}(k_0 = 0)\left[\frac{1}{|W|}\sum_{w \in W}D_k^{(B_m)}(x,wy)\right] = \frac{1}{|W|}<x,x>g(x,y)
\end{equation*}
For the remaining term, note that both $d$ and $f$ are $W$-invariant. Write $\Delta_k = \sum_{i=1}^m T_i^2$, where $T_i$ is the difference Dunkl operator (see \cite{Ros1} p. 5). Then using (see \cite{Ros1} p. 6), one has the derivation formula $T_i (df) = dT_i(f) + fT_i(d)$. It gives that
\begin{equation*}
\Delta_k(df) = d\Delta_k(f) + f\Delta_k(d) + 2\sum_{i=1}^m(T_i(d))(T_i(f)) 
\end{equation*} 
Moreover, $T_i(d) = \partial_i(d)$ and $T_i(f) = \partial_i(f)$ by $W$-invariance. Next we compute : 
\begin{equation*}
\Delta_k(d)(x,y) = 2k_1d(x,y)\sum_{i\neq j}\frac{1}{y_i}\left[\frac{1}{y_i - y_j} + \frac{1}{y_i+y_j}\right] = 4k_1d(x,y)\sum_{i \neq j}\frac{1}{y_i^2-y_j^2} = 0
\end{equation*}
As a result : 
\begin{align*}
\Delta_k(df)(x,y) &= [d\Delta_k(f)](x,y) + 2\sum_{i=1}^m[(\partial_i(d))(\partial_i(f))](x,y) 
\\& = [d\Delta_k(f)](x,y) + 2[d\sum_{i=1}^m\frac{1}{y_i}(\partial_i(f))](x,y)
\\& = d(x,y)[\Delta_k + 2\sum_{i=1}^m\frac{1}{y_i}\partial_i]f(x,y) = d(x,y)\Delta_k^{(B_m)}(k_0 = 1)f(x,y)
\end{align*}
When $k_0 = 1$, $f$ fits the generalized Bessel function in the $B_m$ case $\Rightarrow \Delta_k(df)(x,y) = (1/|W|)<x,x>df(x,y)$. Finally : 
\begin{equation*}
\Delta_k\left[\frac{1}{|W|}\sum_{w \in W} D_k(x,wy)\right] = <x,x> \left[\frac{1}{|W|}\sum_{w \in W} D_k(x,wy)\right] \qquad \qquad \blacksquare
 \end{equation*}

\section{Alcove-valued process}  
\subsection{$\beta$-Jacobi processes}
Recall that the eigenvalues of the real Jacobi matrix process of parameters $(p,q)$ (see \cite{Dou} for facts on this process) satisfy :
{\small \begin{equation*}
d\lambda_i(t) = 2\sqrt{(\lambda_i(t)(1-\lambda_i(t))}d\nu_i(t) + \left[(p - (p+q)\lambda_i(t)) + 
\sum_{j \neq i}\frac{\lambda_i(t)(1-\lambda_j(t)) + \lambda_j(t)(1-\lambda_i(t))}{\lambda_i(t) - \lambda_j(t)}\right] dt
\end{equation*}}
for $0 < \lambda_m(0) < \dots < \lambda_1(0) < 1$ and all $t  < \inf\{s > 0, \lambda_m(s) = 0 \, \textrm{or}\, \lambda_1(s) =1\} \wedge \tau$.
The $\beta$-Jacobi process is defined as a solution, whenever it exists, of the SDE differring from the one above by a parameter $\beta > 0$ in front of the bracket. It is easy to see that if 
$\lambda$ is a $\beta$-Jacobi process of parameters $(p,q)$, then $1-\lambda$ is a $\beta$-Jacobi process of parameters $(q,p)$. As mentioned in the introductory part, the connection with root systems is not new in its own (\cite{Bee}) however we prefer giving some details of this transition. 
Setting $\lambda_i = \sin^2\phi_i$ then $\phi_i = \arcsin\sqrt{\lambda_i} := s(\lambda_i)$ and $0 < \phi_m < \dots < \phi_1 < \pi/2$. The first and second derivatives of $s$ are given by:
\begin{equation*}
s'(\lambda_i) = \frac{1}{\sin 2\phi_i},\qquad s''(\lambda_i) = \frac{2(2\sin^2\phi_i - 1)}{\sin^3 2\phi_i} = -\frac{2\cos 2\phi_i}{\sin^3 2\phi_i}
\end{equation*} 
Using
\begin{equation*}
\sin^2\phi_i - \sin^2 \phi_j = 2\sin(\phi_i + \phi_j)  \sin(\phi_i - \phi_j)
\end{equation*}
\begin{equation*}
\sin^2\phi_i\cos^2\phi_j + \cos^2\phi_i\sin^2\phi_j = \sin^2(\phi_i+\phi_j) + \sin^2(\phi_i-\phi_j)
\end{equation*}
then, Ito's formula gives: 
\begin{align*}
d\phi_i(t) & = d\nu_i(t) + \beta\frac{(p - (p+q)\sin^2\phi_i)}{\sin 2\phi_i} - \cot 2\phi_i dt 
\\&  + \frac{\beta}{2}\frac{dt}{\sin2\phi_i(t)}\sum_{j \neq i}\frac{\sin^2(\phi_i(t)+\phi_j(t)) + \sin^2(\phi_i(t)-\phi_j(t))}{\sin(\phi_i(t) + \phi_j(t)) \sin(\phi_i(t) - \phi_j(t))}
\end{align*}
Writing $\sin^2\phi_i = (1-\cos 2\phi_i)/2$, then
\begin{align*}
d\phi_i(t) &= d\nu_i(t) + \beta\frac{(p-q)}{2}\frac{dt}{\sin2\phi_i(t)} + \frac{\beta(p+q)-2}{2}\cot 2\phi_i(t)dt 
\\&  + \frac{\beta}{2}\frac{dt}{\sin2\phi_i(t)}\sum_{j \neq i}\frac{\sin^2(\phi_i(t)+\phi_j(t)) + \sin^2(\phi_i(t)-\phi_j(t))}{\sin(\phi_i(t) + \phi_j(t)) \sin(\phi_i(t) - \phi_j(t))}
\end{align*}
where $0< \phi_m(0) < \dots <\phi_1(0) < \pi/2$. Moreover,
\begin{equation*}
\sin2\phi_i = [\cot(\phi_i+\phi_j) + \cot(\phi_i-\phi_j)]\sin(\phi_i+\phi_j)\sin(\phi_i-\phi_j)
\end{equation*} 
which gives 
\begin{align*}
d\phi_i(t) & = d\nu_i(t) + \beta\frac{(p-q)}{2}\frac{dt}{\sin2\phi_i(t)} + \frac{\beta(p+q)-2}{2}\cot 2\phi_i(t)dt 
\\&  + \frac{\beta}{2}\sum_{j \neq i}\frac{[1/\sin^2(\phi_i(t)+\phi_j(t))] + [1/\sin^2(\phi_i(t)-\phi_j(t))]}{\cot(\phi_i(t)+\phi_j(t)) + \cot(\phi_i(t)-\phi_j(t))}dt
\end{align*}
Using $1+\cot^2z = 1/\sin^2z$, then 
\begin{align*}
&d\phi_i(t) = d\nu_i(t) + \beta\frac{(p-q)}{2}\frac{dt}{\sin2\phi_i(t)} + \frac{\beta(p+q)-2}{2}\cot 2\phi_i(t)dt 
\\&+ \frac{\beta}{2}\sum_{j \neq i}\frac{\cot^2(\phi_i(t)+\phi_j(t)) + \cot^2(\phi_i(t)+\phi_j(t)) + 2}{\cot(\phi_i(t)+\phi_j(t)) + \cot(\phi_i(t)-\phi_j(t))}dt 
\\& =  d\nu_i(t) +  \beta\frac{(p-q)}{2}\frac{dt}{\sin2\phi_i(t)} + \frac{\beta(p+q)-2}{2}\cot 2\phi_i(t)dt + \beta \times 
\\& \sum_{i\neq j}\left\{\frac{1-\cot(\phi_i(t)+\phi_j(t))\cot(\phi_i(t)-\phi_j(t))}{\cot(\phi_i(t)+\phi_j(t)) + \cot(\phi_i(t)-\phi_j(t))}+\frac{\cot(\phi_i(t)+\phi_j(t)) + \cot(\phi_i(t) - \phi_j(t))}{2}\right\}dt 
\end{align*}
Using
\begin{equation*}
-\cot(u+v) = \frac{1-\cot(u) \cot(v)}{\cot(u) + \cot(v)}.
\end{equation*}
and
\begin{equation*}
\frac{1}{\sin 2\phi_i} = \frac{2\cos^2\phi_i - \cos 2\phi_i}{2\sin\phi_i\cos\phi_i} = \cot \phi_i - \cot2\phi_i 
\end{equation*}
we finally obtain 
\begin{align}\label{JE}
d\phi_i(t) = d\nu_i(t)  +\left[k_0\cot \phi_i + k_1 \cot 2\phi_i(t)dt +  k_2\sum_{i\neq j}[\cot(\phi_i+\phi_j) + \cot(\phi_i -\phi_j)]\right]dt
\end{align}where 
\begin{equation}\label{MF}
2k_0 = \beta(p-q),\quad k_1 = \beta(q-(m-1)) -1, \quad 2k_2 = \beta.
\end{equation}
Easy computations show that $\pi/2 - \phi$ satisfies (\ref{JE}) with $(p,q)$ intertwined.

\subsection{Eigenfunctions and Heckman-Opdam's functions.}
Let $k_2 > 0$ and $\mathscr{L}$ be the generator of $\phi$, then the eigenfunctions of $\mathscr{L}$ are given by Gauss hypergeometric series: in fact, let $\mathscr{A}$ be the generator of $(\lambda_1,\dots,\lambda_m)$ 
(see \cite{Dou} p. 135 for $\beta=1$) : 
\begin{align*}
&\mathscr{A}  = 2\sum_{i=1}^m\lambda_i(1- \lambda_i)\partial_i^2 + \beta \sum_{i=1}^m \left[p-(p+q)\lambda_i + \sum_{j \neq i}\frac{\lambda_i(1-\lambda_j) + \lambda_j(1-\lambda_i)}{\lambda_i - \lambda_j}\right]\partial_i
\\& = 2\sum_{i=1}^m\lambda_i(1- \lambda_i)\partial_i^2 + \beta \sum_{i=1}^m[p - (m-1)-(p+q - 2(m-1))\lambda_i]\partial_i + 2\beta\sum_{i \neq j}\frac{\lambda_i(1-\lambda_i)}{\lambda_i - \lambda_j}\partial_i
\end{align*}
From Equation (\ref{GF}) ($k_2$ plays the role of $k_1$), one can see that ${}_2F_1^{(1/k_2)}(a,b,c;\lambda)$ is the unique symmetric analytic function $u$ such that $u(0) = 1$ which satisfies 
\begin{equation*}
\mathscr{A}u(\lambda) = 2mab\,u(\lambda), \quad  2c = \beta p = 2k_0 + k_1 + 2k_2(m-1) +1, \,2a+2b +1 -2c = k_1.
\end{equation*}  
with $k_i,\,0 \leq i \leq 2$ cited in (\ref{MF}). Setting $\sin^2\phi := (\sin^2\phi_1,\dots,\sin^2\phi_m)$, then $\mathscr{A}$ transforms to $\mathscr{L}$. Hence: 
\begin{equation*}
\mathscr{L}[u(\sin^2\phi)] = 2mab\,[u(\sin^2\phi)] 
\end{equation*} 
In the same spirit, one can also interpret $\mathscr{L}$ as the ``radial part'' of the trigonometric version Dunkl-Cherednik Laplacian (with $\cot$ replacing $\coth$, 
\cite{Bee},\cite{Opdam}). By radial part, we mean the restriction on $W$-invariant functions. Besides, this Laplacian arises, as for Dunkl and Cherednik-Dunkl ones, from differential-difference first-order operators. However, this comes beyond the spirit of this work and will not be done here.   

\subsection{Existence and uniqueness of a strong solution.}
The involved root system is the non reduced $BC_m$ defined by 
\begin{eqnarray*}
R &=& \{\pm e_i, \, \pm 2e_i,\, 1 \leq i \leq m,\, \pm(e_i \pm e_j),\, 1 \leq i < j \leq m\} \\
R_+ & = &  \{e_i, \,  2e_i,\, 1 \leq i \leq m,\, (e_i \pm e_j),\, 1 \leq i < j \leq m\} \\
\Delta &=& \{e_i - e_{i+1}, \,1 \leq i \leq m-1,\, e_m\}
\end{eqnarray*}  
When $k_0 = 0 (p=q)$, it reduces to the reduced $C_m$ system  
\begin{eqnarray*}
R&=& \{\pm e_i \pm e_j,\,1 \leq i < j \leq m,\, \pm 2e_i, \,1 \leq i \leq m\}\\
R_+&=& \{e_i \pm e_j,\,1 \leq i < j \leq m,\, 2e_i, \,1 \leq i \leq m\}\\
\Delta &=& \{e_i - e_{i+1}, \,1 \leq i \leq m-1,\, 2e_m\}
\end{eqnarray*}
and it is known as the ultraspheric case. The Weyl group action on $\re^m$ gives rise to three orbits so that the multiplicity function is given by $k = (k_0,k_1,k_2)$. Setting 
$\tilde{\phi}_i := \phi_i/\pi$, then the process is valued in the {\it positive Weyl alcove} (see \cite{Hum}) defined by : 
\begin{equation*}
\tilde{A} = \{\tilde{\phi} \in \re^m, \, < \alpha,\tilde{\phi}> \, > 0 \, \forall \alpha \in \Delta\, <\tilde{\alpha},\tilde{\phi}>\,  < 1\}
\end{equation*}
where $\tilde{\alpha} = 2e_1$ is the highest positive root (that is $\tilde{\alpha} - \alpha \in R_+ \, \forall \alpha \in R$, see \cite{Hum}). 
The associated affine Weyl group $W_a$ is the semi-direct product of $W$ and the translation group corresponding to the coroot lattice 
($\mathbb{Z}$-span of $\{2\alpha/||\alpha||^2,\, \alpha \in R\})$. 
The generator writes in this case: 
\begin{equation*}
\mathscr{L}g(\phi) := \frac{1}{2}\Delta g(\phi) - <\nabla g(\phi), \nabla\Phi(\phi)>,\quad \Phi(\phi) = -\sum_{\alpha \in R_+}k(\alpha)\log\sin(<\alpha,\phi>)
\end{equation*}
Thus, with minor modifications, Theorem \ref{exi} states that (\ref{JE}) has a unique strong solution for all $t > 0$ subject to $k_0 >0,\,k_1 > 0,\, k_2 > 0 \Leftrightarrow \beta > 0,\, 
p > q > (m-1) + 1/\beta$. Applying this to $\pi/2-\phi$, this holds for $\beta > 0,\, q > p > (m-1) + 1/\beta$. Since the ultraspheric case still involves a root system, then (\ref{JE}) has a unique strong solution for $p \wedge q > (m-1) + 1/\beta$ which simplifies to $p \wedge q > m$ in the real case $\beta = 1$ and $p \wedge q > m-1/2$ in the complex one $\beta =2$. 
Theorem \ref{exi} is modified as follows: $\partial \tilde{A} = \cup_{\alpha \in \Delta}H_{\alpha} \cup H_{\alpha,1}$ where 
\begin{equation*}
H_{\alpha,1} = \{\tilde{\phi}, \, <\tilde{\alpha},\tilde{\phi}> = 1\} =  \{\phi, \, \pi - <\tilde{\alpha},\phi> = 0\}
\end{equation*}
Compared with (\ref{DE1}), the convex function $x \mapsto -\ln(<\alpha,x>)$ should be substituted by $\phi \mapsto -\ln(\sin(<\alpha,\phi>))$ and one has to deal with an additional term in the expression of the boundary process $(L_t)_{t \geq 0}: \displaystyle {\bf 1}_{\{\pi - <\tilde{\alpha},\phi> = 0\}}$. Then the occupation density formula writes:
\begin{align*}
\int_0^{\pi/2}L_t^a(\pi - <\tilde{\alpha},\phi>)|\theta^{'}(a)|da &= <\tilde{\alpha},\tilde{\alpha}>  \int_0^t |\theta^{'}(\pi - <\tilde{\alpha},X_s>)| ds 
\\& = <\tilde{\alpha},\tilde{\alpha}>  \int_0^t |\theta^{'}(<\tilde{\alpha},X_s>)| ds 
\end{align*}
since $\cot(\pi-z) = -\cot(z)$. Hence, the same proof applies and Lemma 1 remains valid for $\tilde{\alpha} \in R_+$. Besides, either it will exist $\alpha \in \Delta$ such that $<\alpha,x> = 0$ and Lemma \ref{Abass} applies, or we will need to prove that $<n(x),\tilde{\alpha}> \neq 0$ if $x$ belongs only to $H_{\tilde{\alpha},1}$. Let us first recall that the highest root is the unique positive root such that $\tilde{\alpha} - \alpha \in R_+$ for all $\alpha \in R_+$. Thus it may be written as  $\tilde{\alpha} = \sum_{\alpha \in \Delta}a_{\alpha}\alpha$ where $a_{\alpha} \geq 1$. Else, if there exists $\alpha_0 \in \Delta$ such that $a_{\alpha_0} < 1$ and since $\tilde{\alpha}$ must be greater than all simple roots (in particular greater than $\alpha_0$) then
\begin{equation*}
\tilde{\alpha} - \alpha_0 = (a_{\alpha_0} - 1)\alpha_0 + \sum_{\alpha_0 \neq \alpha \in \Delta}a_{\alpha}\alpha = c_{\alpha_0}\alpha_0 + \sum_{\alpha_0\neq \alpha \in \Delta} c_{\alpha}\alpha
\end{equation*} 
for some $c_{\alpha} \geq 0$. Our claim follows from the fact that $\Delta$ is a basis. 
Next, it is not difficult to see from the definition of $n(x)$ and the fact that $<\alpha ,x > > 0$ for all $\alpha \in \Delta$ that $n(x)$ is colinear to $-\sum_{\alpha \in \Delta}\alpha$. It follows that 
\begin{equation*}
<n(x),\tilde{\alpha}> = - c\sum_{\alpha \in \Delta} <\alpha,\tilde{\alpha}> = -c \sum_{\alpha \in \Delta}\sum_{\theta \in \Delta}a_{\alpha}<\alpha,\theta> 
\end{equation*}
If $<n(x),\tilde{\alpha}> = 0$, then 
\begin{equation*}
||\sum_{\alpha \in \Delta}\alpha||^2 = \sum_{\alpha \in \Delta}\sum_{\theta \in \Delta}<\alpha,\theta>  \leq \sum_{\alpha \in \Delta}\sum_{\theta \in \Delta}a_{\alpha}<\alpha,\theta>  = 0 
\end{equation*}
which implies that $n(x) = 0$.$\hfill \blacksquare$

\subsection{Brownian motion in Weyl alcoves}
Let
\begin{equation*}
h_1(\phi) := \prod_{\alpha \in R_+}\sin(<\alpha,\phi>)
\end{equation*}
Then, $h_1$ is strictly positive on $\tilde{A}$ and vanishes for $\phi \in \partial A$. One can also show that $(1/2)\Delta h_1 = ch$ for some strictly negative constant $c$. 
Let $P_t^{h_1}$ denote the semi group given by : 
\begin{equation*}
P_t^{h_1}f(\phi) : = e^{-ct} \frac{P_t(h_1f)(\phi)}{h_1(\phi)}, 
\end{equation*} where $P_t$ denotes the semi group of the process consisting of $m$-independent BMs in $A$ killed when it first reaches $\partial A$. The corresponding generator writes : 
\begin{equation*}
\mathscr{L}^{h_1}(f) = \frac{1}{h_1}\left[\frac{1}{2}\Delta - c\right](h_1f) = \frac{1}{2}\Delta + \sum_{i=1}^m(\partial_i\log h_1)\partial_i f
\end{equation*}
which fits our generator for $k_2 = 1\,(\beta = 2), k_1 = 2\,(q = m+1/2),\, k_0 = 1 (p = q+1 = m+3/2)$. In the ultraspheric case, this becomes $\beta = 2,\,p=q = m+1/2$. In both cases, these parameters correspond to the process consisting of $m$ BMs constrained to stay in the $BC_m$ and $C_m$- Weyl alcoves respectively. Note that $p,q$ are not integers which means that these processes BM can not be realized as eigenvalues processes of complex matrix Jacobi processes which is also the case for the BM in the $B_m$-Weyl chamber since 
$\delta = m+1/2$. 

\subsection{The first hitting time $\tilde{T}_0$}
We define similarly the first hitting time of alcove's walls by $\tilde{T}_0 = \inf\{t > 0,\, (\phi(t)/\pi) \in \partial \tilde{A}\} = 
\tilde{T}_{\tilde{\alpha}} \wedge \inf \{\tilde{T}_{\alpha},\,\alpha \in \Delta\}$, where 
\begin{eqnarray*}
\tilde{T}_{\alpha} &:=& \inf\{t > 0,\, <\alpha,\phi(t)> = 0\},\\
\tilde{T}_{\tilde{\alpha}} &:= & \inf\{t > 0, <\tilde{\alpha}, \phi(t)> = 2\phi_1 = \pi\}, 
\end{eqnarray*}
and $\phi$ is the unique strong solution for all $t \geq 0$ of \footnote{$k(2e_i) = k_1/2$ for all $1 \leq i \leq m$.}:  
\begin{equation*}
d\phi(t) = d\nu(t) + \sum_{\alpha \in R+}k(\alpha)\cot(<\alpha,\phi(t)>)\alpha \, dt, \quad \frac{\phi(0)}{\pi}   \in \tilde{A},  
\end{equation*}
for the non reduced root system $R=BC_m$ with $k(\alpha) > 0$ for all $\alpha$ and $p \wedge q > (m-1) + 1/\beta$. Let us focus on $\tilde{T}_{\alpha_0}$ for some $\alpha_0 \in \Delta$. We shall distinguish two cases : 

\subsubsection{$\alpha_0 = e_i - e_{i+1}, \, 1 \leq i \leq m-1$}
The same scheme described in the proof of Proposition \ref{Sacha} applies here since the main ingredients used there are the SDE and the fact that $\sigma_0(\alpha) \in R_+$ if 
$\alpha \neq \alpha_0$. The second assertion follows from $\sigma_0(2e_j) = 2\sigma_0(e_j) = 2(\delta_{ij}e_{i+1} + \delta_{(i+1)j}e_{i} + {\bf 1}_{\{j\neq i,j \neq i+1\}}e_j)  \in R_+$. 
As a result, one writes for all $t \geq 0$: 
\begin{align*}
d<\alpha_0,\phi(t)> = ||\alpha_0||d\gamma_t + k_2||\alpha_0||^2\cot<\alpha_0,\phi(t)> dt+ \sum_{\substack{\alpha \in R_+\\ \alpha \neq \alpha_0}}k(\alpha)a(\alpha)\cot<\alpha,\phi(t)>dt
\end{align*} 
where $a(\alpha) = <\alpha_0,\alpha>$. 
\begin{align*}
d<\alpha_0,\phi(t)> = ||\alpha_0||d\gamma_t + k_2||\alpha_0||^2\cot(<\alpha_0,\phi(t)>) dt + F_t
\end{align*}  where 
\begin{align*}
F_t = \sum_{\substack{\alpha \in R_+\setminus \alpha_0 \\  a(\alpha) >0}}k(\alpha)a(\alpha)[\cot(<\alpha,\phi(t)>) - \cot(<\sigma_0(\alpha),\phi(t)>)],
\end{align*}
where $\sigma_0 = \sigma_{\alpha_0}$. This drift is strictly negative on $\{\tilde{T}_{\alpha_0} = \infty\}$ since $\phi \mapsto \cot \phi$ is a decreasing function, $<\alpha_0, \phi(t)> > 0$ and since : 
\begin{equation*}
<\alpha - \sigma_0(\alpha),\phi(t)> = 2\frac{a(\alpha)}{||\alpha||^2}<\alpha_0,\phi(t)> > 0 .
\end{equation*}
This implies that $\mathbb{P}_{x}(\forall t \geq 0,\, <\alpha_0,\phi(t)> \leq Z_{t}) = 1$ where $\phi(0) = x$ and : 
\begin{equation*}
dZ_t = ||\alpha_0||d\gamma_t + ||\alpha_0||^2k_2\cot(Z_t) dt, \quad Z_0 = <\alpha_0,\phi(0)> = x
\end{equation*}
on the same probability space. Using (\ref{JE}) with $\beta = 1, m=1$, one can easily see that $(Z_{t})_{t \geq 0} = (\arcsin\sqrt{J}_{||\alpha_0||^2t})_{t \geq 0}$ where $J$ is a one dimensional Jacobi process of parameters $d = 2k_2+1, d' = 1$ (see \cite{War}) :   
that is : 
\begin{equation*}
dJ_t = 2\sqrt{J_t(1-J_t)}d\gamma_t + (d - (d+1)J_t) dt,  \quad 0 < k_2 < 1/2 \Leftrightarrow 0 < d < 2. 
\end{equation*}  
As $J$ hits $0$ almost surely when $0 < d < 2$ (use the skew product in \cite{War} and properties of squared Bessel processes), then so does $Z$ and by the way $<\alpha_0,\phi>$ for 
$k_2 < 1/2 \Rightarrow \tilde{T}_{\alpha_0} < \infty$ a. s. 
\subsubsection{$\alpha_0 = e_m$} 
Compared with the previous case, the difference arises from the fact that $\sigma_0(\alpha) \in R_+$ if $\alpha \in R_+ \setminus \{e_m,2e_m\}$ and the latter is easily checked since for $\alpha = e_i \pm e_j$ this amounts to consider the reduced root system $B_m$, else for $\alpha = e_i$ with $i \neq m, \,\sigma_0(e_i) = e_i$. According to this, one gets : 
\begin{align*}
d<\alpha_0,\phi(t)> = d\phi_m(t) = d\gamma_t + k_0\cot(\phi_m(t)) dt + k_1\cot(2\phi_m)+F_t
\end{align*}  where 
\begin{align*}
F_t &= \sum_{\substack{\alpha \in R_+\setminus \{e_m,2e_m\} \\  a(\alpha) >0}}k(\alpha)a(\alpha)[\cot(<\alpha,\phi(t)>) - \cot(<\sigma_0(\alpha),\phi(t)>)]
\end{align*}
where $R_+^1 = \{e_i - e_j,\, 1 \leq i < j \leq m\}$. Using once again (\ref{JE}), we shall compare this process with $(\arcsin\sqrt{J_t})_{t \geq 0}$ where 
\begin{equation*}
dJ_t = 2\sqrt{J_t(1-J_t)}d\gamma_t + (d - (d+d')J_t) dt,  \quad d'= k_1 + 1,\, d = 2k_0 + k_1 + 1. 
\end{equation*}
Hence, $\tilde{T}_{e_m} < \infty$ a.s. if $0< 2k_0 + k_1< 1/2 \Leftrightarrow \beta p - (\beta(m-1) < 2 $. This agrees with the case $m=1$ for which $p < 2$ (use the skew product in \cite{War}). Finally, note that since $a(\alpha) = 0$ for $\alpha \in \{e_i, 2e_i,\, 1 \leq i \leq m-1\}$, $F$ only involves $k_2 = \beta$ which is independent from $p,q$. Keeping in mind that 
$\pi/2 - \phi$ is still a $\beta$-Jacobi process with $(p,q)$ intertwined which has no effect on the strict negativity of $F$ by the above remark, we conclude that 
$\tilde{T}_{\tilde{\alpha}} < \infty$ for $0 < \beta q - \beta(m-1) < 2$.
$\hfill \blacksquare$

\subsection{Semi-group density}
We end this paper by giving the semi group density of the $\beta$-Jacobi process. Before proceeding, we briefly consider two cases for which we can write down the semi-group density: the univariate case and the complex Hermitian one $(\beta=2)$. Let $\jac$ denote the Jacobi polynomial of degree $n$ defined by (\cite{Askey}): 
\begin{equation*}
\jac(\lambda) := \frac{(r+1)_n}{n!}{}_2F_1\left(-n,n+r+s+1,r+1; \frac{1-\lambda}{2}\right), 
\end{equation*}
for $\lambda \in [-1,1], \,r,s > -1$, where ${}_2F_1$ is the univariate Gauss hypergeometric function. 
These polynomials are orthogonal with respect to the measure $Z^{r,s}(\lambda)d\lambda := (1-\lambda)^{r}(1+\lambda)^{s}d\lambda$ and the associated inner product in $L^2([-1,1])$ given by 
\begin{equation*}
<f,g>_{L^2([-1,1])} : = \int_{[-1,1]}f(\lambda)g(\lambda)Z^{r,s}(\lambda)d\lambda
\end{equation*}
Moreover, $(\jac)_{n \geq 0}$ form a complete set of this Hilbert space and satisfy 
\begin{equation*}
\left\{\sqrt{1-\lambda^2}\partial_\lambda^2 + [(s - r) - (s +r +2)\lambda]\partial_\lambda\right\}\jac(\lambda) = -n(n+r+s+1)\jac(\lambda)
\end{equation*}   
The above eigenoperator defines a diffusion which is related to the one we considered with $m=1$ via the map $\lambda \mapsto (1-\lambda)/2$ and a deterministic time change ($t \mapsto t/2)$. 
The semi group density w.r.t Lebesgue measure is written (see \cite{Wong},\cite{Dem1})
\begin{equation*}
p_t^{r,s}(\theta,\lambda) = \sum_{n=0}^{\infty}e^{-2r_nt}\jac(\theta)\jac(\lambda) W^{r,s}(\lambda)
\end{equation*} 
where $r_n$ denotes the eigenvalues above, $(\jac)_n$ are orthonormal polynomials, $p = 2(r+1),q= 2(s+1)$ and $W^{r,s}(\lambda)d\lambda$ is the probability measure corresponding to the measure 
$Z^{r,s}(\lambda)d\lambda$. No closed forms seems to be known for this density, nonetheless an attempt to get a handier expression was tried in \cite{Dem1}. Multivariate analogs appeared in literature (\cite{Bee}, \cite{Kan}, \cite{Lass1} for instance) and are obtained by applying the Gram-Schmidt orthogonalization to the symmetric Jack polynomials w.r.t. measure
\begin{equation*}
Z_m^{r,s,\beta}(\lambda)d\lambda := \prod_{i=1}^m\lambda_i^{r}(1-\lambda_i)^{s} \prod_{1\leq i < j\leq m}|\lambda_i - \lambda_j|^{\beta}d\lambda_1\dots d\lambda_m
\end{equation*}   
We shall denote them\footnote{The normalization is different from the one used in both \cite{Bee} and \cite{Lass1}.} by $\jam$ for a given partition $\tau$ (instead of $G_{\tau}^{\alpha,\beta}$ used in literature) and stress that some of the properties cited above extend to the higher dimensional case (\cite{Lass1}): an expansion in terms of ${}_2F_1^{(2/\beta)}(-l,b,c,\lambda)$ exists for $\tau = (l^m)$ with $m$ components all equal to $l$; $(\jam)$, where $\tau$ is a partition of length $\leq m$, form a basis of the Hilbert space $L^2([0,1]^m, W_m^{r,s,\beta}(\lambda) d\lambda)$ where $W_m^{r,s,\beta}(\lambda) d\lambda$ is the normalization $Z_m^{r,s,\beta}(\lambda)d\lambda$ in order to be a probability measure (\cite{Lass1}).
The normalizing constant is given by a McDonald-Selberg integral computed in \cite{Kan}. Moreover, $(\jam)_{\tau}$ are the unique symmetric polynomial eigenfunctions of the Laplace Beltrami operator 
$-\mathscr{L}$ (thus defined on $[0,1]^m$) with $\beta (p - (m-1)) = 2(r+1),\,\beta(q-(m-1)) = 2(s+1)$, associated with the eigenvalues   
\begin{equation}\label{vp}
2r_{n,\tau}^{\beta} := 2\left[\sum_{i=1}^m\tau_i(\tau_i-1 -\beta(i-1)) + |\tau|(r+s+\beta(m-1) +2)\right],\,|\tau|=n.
\end{equation}
However, with regard to the strong uniqueness for all $t \geq 0$ previously derived, we shall restrict ourselves to $p \wedge q > (m-1) +1/\beta$. $\beta(q - (m-1)) > 1$ is equivalent to $s > -1/2$ and $\beta(p - (m-1)) > 1$ is equivalent to $r > -1/2$. As a result, $r,s > -1/2$.\\
It is known that the eigenvalues process of the complex Hermitian Jacobi process (or $2$-Jacobi process) is the $h$-transform (in the Doob sense) for $h=V$ of a process whose components are real Jacobi processes of parameters $2(p-(m-1)) = 2(r+1), 2(q-(m-1)) = 2(s+1)$ constrained to never collide (or to stay in the $A_{m-1}$-type Weyl chamber). Here, $V$ denotes as usual the Vandermonde function. More precisely, $V$ is an eigenfunction of the generator of the one dimensional Jacobi process of parameters $(p,q)$ (see appendix in \cite{Dou}), say 
$\mathbf{L}$, that is
\begin{equation*}
\mathbf{L}V = cV = -m(m-1)\left(\frac{2(m-2)}{3} + \frac{p+q}{2}\right) V
\end{equation*}
Noting that the parameters $r,s$ are the same both in the univariate and in the multivariate cases, it follows by Karlin-McGregor formula (\cite{Kar}) that the semi group density writes 
on $\{0 < \lambda_m < \dots < \lambda_1 < 1\}$
\begin{align*}
&K_t^{r,s,2}(\theta,\lambda) := e^{-ct}\frac{V(\lambda)}{V(\theta)} \det\left(\sum_{n=0}^{\infty}e^{-2n(n+r+s+1)t}\jac(\theta_i)\jac(\lambda_j)W^{r,s}(\lambda_j)\right)_{i,j}
\\& = e^{-ct}\det\left(\sum_{n=0}^{\infty}e^{-2n(n+r+s+1)t}\jac(\theta_i)\jac(\lambda_j)\right)_{i,j}\frac{W_m^{r,s,2}(\lambda)}{V(\theta)V(\lambda)}
\\& = e^{-ct}\left[\sum_{\sigma_1 \in S_m}\epsilon(\sigma_1)\sum_{n_1,\dots,n_m \geq 0}e^{-2\sum_{i=1}^mn_i(n_i + r + s+1)t}\prod_{i=1}^mP_{n_i}^{r,s}(\theta_i)P_{n_{i}}^{r,s}
(\lambda_{\sigma_1(i)})\right]\frac{W_m^{r,s,2}(\lambda)}{V(\theta)V(\lambda)}
\\& = e^{-ct}\left[\sum_{\sigma_1,\sigma_2 \in S_m}\epsilon(\sigma_1)\sum_{n_1 \geq \dots \geq n_m \geq 0}e^{-2\sum_{i=1}^m n_{\sigma_2(i)}(n_{\sigma_2(i)} + r + s+1)t}
\prod_{i=1}^mP_{n_{\sigma_2(i)}}^{r,s}(\theta_i)P_{n_{\sigma_2(i)}}^{r,s}(\lambda_{\sigma_1(i)})\right]\frac{W_m^{r,s,2}(\lambda)}{V(\theta)V(\lambda)}
\end{align*}
Note that, for a given partition $(n_1\geq \dots \geq n_m \geq 0)$ and a permutation $\sigma_2 \in S_m$, one has 
\begin{equation*}
\sum_{i=1}^m n_{\sigma_2(i)}(n_{\sigma_2(i)} + r + s+1) = \sum_{i=1}^mn_i(n_i + r + s+1)
\end{equation*}
Thus summing first over $\sigma_1$ with the change of variables $\sigma = \sigma_1\sigma_2$, one gets: 
\begin{align*}
K_t^{r,s,2}(\theta,\lambda) &= e^{-ct}\sum_{n_1 \geq \dots \geq n_m \geq 0}e^{-2\sum_{i=1}^mn_{i}(n_{i} + r + s+1)t} 
\frac{\det[P_{n_i}^{r,s}(\theta_j)]_{i,j}}{V(\theta)}\frac{\det[P_{n_i}^{r,s}(\lambda_j)]_{i,j}}{V(\lambda)}W_m^{r,s,2}(\lambda)
\\& = e^{-ct}\sum_{n_1 > \dots > n_m \geq 0}e^{-2\sum_{i=1}^mn_{i}(n_{i} + r + s+1)t} 
\frac{\det[P_{n_i}^{r,s}(\theta_j)]_{i,j}}{V(\theta)}\frac{\det[P_{n_i}^{r,s}(\lambda_j)]_{i,j}}{V(\lambda)}W_m^{r,s,2}(\lambda)
\end{align*} 
Set $n_i = \tau_i + m - i$, then $\tau_1 > \geq \dots \geq \tau_m \geq 0$. Moreover, with regard to (\ref{vp}), one easily check that 
\begin{equation*}
r_{n,\tau}^2 = \sum_{i=1}^m\tau_i(\tau_i + r+s+1 + 2(m-i))
\end{equation*}
so that 
\begin{align*}
\sum_{i=1}^mn_{i}(n_{i} + r + s+1) = r_{n,\tau}^2 - c/2
\end{align*}
The final result writes 
\begin{align*}
K_t^{r,s,2}(\theta,\lambda) &= \sum_{\tau_1 \geq \dots \tau_m \geq 0}e^{-2r_{n,\tau}^2} 
\frac{\det[P_{\tau_i+ m-i}^{r,s}(\theta_j)]_{i,j}}{V(\theta)}\frac{\det[P_{\tau_i+m-i}^{r,s}(\lambda_j)]_{i,j}}{V(\lambda)}W_m^{r,s,2}(\lambda)
\\& = \sum_{\tau_1 \geq \dots \tau_m \geq 0}e^{-2r_{n,\tau}^2} P_{\tau}^{r,s,2}(\theta)P_{\tau}^{r,s,2}(\lambda)W_m^{r,s,2}(\lambda)
\end{align*}
where we used the determinantal representation of the Jacobi multivariate polynomials in the complex case{\footnote{We adopt a different normalization since we consider orthonormal polynomials.} (see \cite{Lass1}) : 
\begin{equation*}
P_{\tau}^{r,s,2}(\lambda) = \frac{\det[P_{\tau_i+m-i}^{r,s}(\lambda_j)]_{i,j}}{V(\lambda)}
\end{equation*}

From these observations, it is natural to claim that : 
\begin{pro}
The semi group density of the $\beta$-Jacobi process is given by 
\begin{equation}
\label{densite}
K_t^{r,s,\beta}(\theta, \lambda) := \sum_{n=0}^{\infty}\sum_{|\tau|=n}e^{-r_{n,\tau}t}\jam(\theta)\jam(\lambda) W_m^{r,s}(\lambda){\bf 1}_{\{0 < \lambda_m < \dots < \lambda_1 < 1\}}
\end{equation}
with respect to $d\lambda$. As a result, it is positive.
\end{pro}
{\it Proof}: given a bounded symmetric function $f$ on $[0,1]^m$, define
\begin{align*}
T_tf(\theta) &:= \int_{0 < \lambda_m< \dots < \lambda_1 < 1}f(\lambda)\sum_{n=0}^{\infty}\sum_{|\tau|=n}e^{-2r_{n,\tau}^{\beta}t}\jam(\theta)\jam(\lambda) W_m^{r,s}(\lambda)d\lambda 
\end{align*}
for $\theta = (0 < \theta_1<\dots < \theta_m < 1)$ and $T_0f = f$. The above expression makes sense: this uses the boundness of $f$, the exponential term with strictly positive $t$ and Fubini Theorem. Besides, $T_t{\bf 1} = 1$ and $||T_t||$ is bounded for all $t \geq 0$. The first claim follows easily from the orthogonality of the Jacobi polynomials and $P_0 = {\bf 1}$ so that the only non zero term is that correponding to $n=0$. The second one is obvious for $t=0$ and uses the exponential term when $t \geq \epsilon >0$. 
One also easily checks that $T_tT_s = T_{t+s}$ and that $\mathscr{L}T_tf(\lambda) = \partial_tT_tf(\lambda)$ using the dominated convergence theorem. 
Now, let us consider the Cauchy problem associated to $\mathscr{L}$ : 
\begin{equation*}
\left\{\begin{array}{l}
 \displaystyle \frac{\partial u_{f}}{\partial t}(t,\lambda) = \mathscr{L}u_{f}(t,\lambda) \\ 
 u_{f}(0,\cdot) = f, 
\end{array}\right.
\end{equation*}
where $u_{f} \in C^{1,2}(\re_+^{\star} \times S := \{0 < \lambda _m < \dots < \lambda_1 < 1\}) \cap C_b(\re^+ \cap S)$ with reflecting boundary condition : 
\begin{equation*}
< \nabla u(t,\lambda), n(\lambda)> = 0 ,\quad (t,\lambda) \in \re_+^{\star} \times \partial S 
\end{equation*}
where $n(\lambda)$ is a unitary inward normal vector at Define $u_t (f)(\lambda) := u_{f}(t,\lambda)$. It is shown (\cite{Str}) that there is a unique solution to the Cauchy problem with initial condition. Consequently, $(T_t)_{t \geq 0}$ is the semi group of the eigenvalues process $(\lambda(t))_{t \geq 0}$ with density given by $K_t^{r,s,\beta}$. $\hfill \blacksquare$ 

\begin{nota}
As the reader can check, the computations performed in the complex Hermitian case do not restrict to Jacobi polynomials. We only used the determinantal representation in terms of their univariate counterparts. As a result, one gets similar formulas replacing Jacobi by Hermite and Laguerre polynomials.  
\end{nota}
Now, we are able to answer some open questions left in \cite{Dou}. For the real Jacobi matrix ($\beta =1$), it is known that for $p \wedge q \geq m-1$ and if the eigenvalues are distinct at time $t=0$, then they remain distinct forever. It is then natural to wonder if this remains valid when starting from non distinct eigenvalues (see \cite{Dou} p. 138-139). The Markov property together with the previous result for distinct eigenvalues are sufficient to claim that this is true provided that the eigenvalues semi group has a density which is absolutely continuous with respect to Lebesgue measure on $\re^m$. 
By virtue of $K_t^{r,s,1}(\theta,\phi)$, for $p \wedge q > m$, 
\begin{equation*}
\p_{\lambda(0)}(\forall t \geq 0,\, \forall i \neq j,\,\lambda_i(t) \neq \lambda_j(t)) = 1,\quad \lambda_1(0) \geq \dots \geq\lambda_m(0).
\end{equation*}
We argue in the same way to claim that for $p \wedge q \geq m+1$, the process will never hit the boundaries ($0$ and $1$ for $\lambda$ or $0$ and $\pi/2$ for $\phi$) even if it did at time $t=0$.\\
{\bf Acknowledgment} : the author would to thank C. Donati Martin for useful remarks and her careful reading of the paper, and P. Bougerol for explanations of some facts on root systems. A special thank to M. Yor for his intensive reading of the manuscript.

\end{document}